\documentclass[11pt]{article}

\usepackage{amssymb}
\usepackage{amsthm}
\usepackage{amsmath,amsfonts,amssymb}

\allowdisplaybreaks
\usepackage[mathscr]{eucal}
\usepackage{exscale}
\usepackage[square,sort,comma,numbers]{natbib}
\setcitestyle{authoryear,round}

\usepackage{bm}
\usepackage{eqlist}
\usepackage[dvipsnames]{color}
\usepackage[Lenny]{fncychap}
\usepackage{multirow}
\usepackage{graphicx}
\usepackage{algpseudocode}
\usepackage{algorithm}
\usepackage{caption}
\usepackage{hyperref}
\usepackage{framed}
\usepackage{color}
\usepackage{booktabs}
\usepackage{makecell}

\usepackage{framed}
\usepackage{algpseudocode}
\usepackage{indentfirst} 
\usepackage{longtable}

\usepackage{tabu}
\usepackage{tikz}

\usepackage{chngpage}

\usepackage[utf8]{inputenc} 
\usepackage[english]{babel} 
\usepackage[T1]{fontenc}    
\usepackage{amssymb,amsmath,amsfonts} 
\usepackage{braket} 
\usepackage{caption}
\usepackage{algorithm}

\hypersetup{
	colorlinks=true,
	linkcolor=blue,
	citecolor=blue,
	citebordercolor={1 1 0},
}

\newtheorem{theorem}{Theorem}[section]

\newtheorem{definition}[theorem]{Definition}
\newtheorem{remark}[theorem]{Remark}

\font\bigbf=cmbx10 scaled \magstep3

\begin{document}
	
	\title{\bigbf  Trip-based mobile sensor deployment for drive-by sensing with bus fleets}
	
	\author{Wen Ji
		\quad Ke Han$\thanks{Corresponding author, e-mail: kehan@swjtu.edu.cn;}$
		\quad Tao Liu
		\\\\
		\textit{\small Institute of System Science and Engineering, School of Transportation and Logistics,}\\
		\textit{\small Southwest Jiaotong University, Chengdu 611756, China}
	}
	
	\maketitle
	
	\begin{abstract}
	Drive-by sensing (i.e. vehicle-based mobile sensing) is an emerging data collection paradigm that leverages vehicle mobilities to scan a city at low costs. It represents a positive social externality of urban transport activities. Bus transit systems are widely considered in drive-by sensing due to extensive spatial coverage, reliable operations, and low maintenance costs. It is critical for the underlying monitoring scenario (e.g. air quality, traffic state, and road roughness) to assign a limited number of sensors to a bus fleet to ensure their optimal spatial-temporal distribution. In this paper we present a trip-based sensor deployment problem, which explicitly considers timetabled trips that must be executed by the fleet while a portion of them perform sensing tasks. To address the computational challenge in large-scale instances, we design a multi-stage solution framework that decouples the spatial-temporal structures of the sensing task through line pre-selection and bi-level optimization. As a result, the computational complexity is reduced to be sub-linear w.r.t. the number of lines, rather than combinatorial w.r.t. the number of buses in existing vehicle-based approaches. A real-world case study covering 400 km$^2$ in central Chengdu demonstrates the effectiveness of the model in solving large-scale problems. It is found that coordinating bus scheduling and sensing tasks can substantially increase the spatial-temporal sensing coverage. We also provide a few model extensions and recommendation for practice regarding the application of this method. 
	\end{abstract}


	\section{Introduction}\label{secIntro}
	
	Ubiquitous sensing plays a vital role in smart transport or smart city applications. In the past decade, the proliferation of wireless and miniaturized sensing technology has not only made pervasive sensing possible, but also poses considerable challenge in the design, deployment and management of such systems \citep{DSXVF2019}. Recently, drive-by sensing (i.e. vehicle-based mobile sensing) has attracted much attention as an efficient and affordable sensing paradigm \citep{MZY2014,Lee2010}. It leverages the mobility of transport vehicles such as taxis, buses and trucks to scan a city at high resolution \citep{LBDNT2005}. Drive-by sensing has been widely adopted in various urban monitoring scenarios such as air quality \citep{SHS2021,Messier2018}, traffic state \citep{ZLMH2014, GQDRJ2022, YFLRL2021}, noise pollution \citep{Cruz2020}, heat island phenomena \citep{TRMSO2020,FBRDRAP2021}, and infrastructure health \citep{AD2017, Eriksson2008}.
 
	Public transport vehicles such as buses are suitable hosts for drive-by sensing because of their easy-to-predict trajectories, extensive spatial-temporal coverage, and low set-up and maintenance costs \citep{CDMFD2018, Cruz2020, Cruz2020b}. A list of bus drive-by sensing studies can be found in the literature review (Section \ref{subseclrbusds}). The quality of sensing coverage offered by bus fleets is highly dependent on the spatial-temporal distribution of the sensors, and the unique mobility pattern of buses gives rise to the interesting and critical problem of {\it mobile sensor deployment} (MSD). Due to the presence of numerous bus lines and even larger number of buses in a dense urban area, the MSD is a complex combinatoric optimization problem encompassing both spatial and temporal dimensions, constrained by bus operational considerations.

	Existing studies on MSD  for bus fleets can be categorized as {\it line-based} and {\it vehicle-based}. The first kind solely focuses on spatial coverage by designing algorithms to allocate sensors to various bus lines instead of individual buses \citep{AD2017, Kaivonen2020}. The resulting optimization problem is usually tractable, but the coarse granularity of sensor deployment does not promise a good sensing outcome for the lack of consideration of bus movements. The second kind concerns with both spatial and temporal coverage by allocating sensors to individual buses based on their real-world GPS trajectories, using heuristic or approximate algorithms \citep{Gao2016, TRMSO2020}. These approaches, typically formulated as subset selection problems,  suffer from not only computational burden due to the large number of candidate buses, but also uncertainties in their trajectories or services, having no input on how bus fleets are operated.

	This paper addresses (1) the aforementioned shortcomings of line-based and vehicle-based  approaches, and (2) the computational challenge in large-scale instances, by making the following contributions: 
	
\begin{itemize}
\item {\bf Model formulation.} Unlike existing studies where sensors are either allocated to lines or vehicles, we propose a trip-based sensor deployment approach, which explicitly considers timetabled trips that must be executed by the fleet while a portion of them perform sensing tasks. This approach ensures that deployed sensors are optimally circulated in the area without affecting the level of bus service. It also allows several extensions with real-world considerations (Section \ref{subsecdiscussion}).

\item {\bf Solution algorithm.} A multi-stage solution framework, including a sequential and a joint approach, is developed, which encompasses bus line pre-selection, minimum fleet size, and trip chain formation. The sequential approach is treated as a benchmark as it pertains to the subset selection paradigm found in the literature. The joint approach, on the other hand, considerably reduces the computational complexity by decoupling the spatial-temporal structures of the sensing task. As a result, the computational burden grows sub-linearly w.r.t. the number of bus lines, as opposed to the combinatorial nature of existing models in the literature.

\item {\bf Practical insights.} A large-scale case study demonstrates the computational efficiency of the proposed methods. It is also found that (1) the non-uniqueness of the minimum fleet solution can be explored to greatly improve sensing quality by coordinating bus scheduling and sensing tasks; (2) the effect of such coordination is more pronounced for high-frequency sensing scenarios; (3) typically 1-2 sensors are sufficient to saturate a bus line, which leads to a practical method that quickly constructs an approximate solution. More specific findings and insights are presented in Section \ref{sec_conclusion}.
\end{itemize}

	This paper explores the trip-based mobile sensor deployment paradigm that encompasses bus scheduling and sensor assignment, while offering decision support tools on a tactical level that are more operationally robust and computationally tractable than existing approaches.

    The remainder of this paper is organized as follows. Section \ref{secLR} reviews some relevant work. The mobile sensor deployment problem is formally defined in Section \ref{secPS}. Section \ref{sec_MLP_two_mode}  develops several solution procedures and model reduction techniques. A real-world case study of our models and algorithms is presented in Section \ref{secCS}, followed by some discussions and recommendations in Section \ref{sec_conclusion}.

    \section{Related work}\label{secLR}
    Work related to this study is divided into three parts: Drive-by sensing, mobile sensor deployment, and bus scheduling.
	
	\subsection{Drive-by sensing}
	
        With the development of wireless and pervasive sensing technologies, drive-by sensing has become a popular data collection method, which offers good coverage in both space and time \citep{MZY2014, ZWXG2014}. Various types of vehicles have been considered, including taxis \citep{Honicky2008}, buses \citep{Castineira2008}, trash trucks \citep{dADKKR2020}, unmanned aerial vehicles \citep{Lambey2021}, and dedicated vehicles \citep{Messier2018, LTH2021}.  Taxis are widely used in urban drive-by sensing, but their spatial-temporal coverage may be biased which mainly be due to factors such as travel demand and revenue \citep{OASSR2019, ZMLL2015}. Public service vehicles such as trash trucks can cover a larger portion of street segments, but at a much lower time resolution, and usually operate only a few days a week \citep{ADRMdR2018}. Unmanned aerial vehicles and dedicated vehicles (cars designed to perform urban sensing) have also been considered for some targeted sensing tasks because of their high flexibility and controllability, but their procurement, deployment and maintenance costs are high \citep{Rashid2020, Messier2018}.

       Bus-based urban sensing has the following advantages: (1) reliable operation due to regular timetables and fixed routes; (2) wide spatial coverage, extending to areas with sparse taxi counts; and (3) low operation and maintenance costs.

	\subsection{Sensor deployment in bus drive-by sensing}\label{subseclrbusds}

        The drive-by sensing capability of bus fleets have been investigated with empirical evidence. \cite{CDMFD2018} analyze the spatial-temporal coverage of bus fleet in Rio de Janeiro, and found that only 18\% of the fleet contribute to 94\% of the total coverage. \cite{Cruz2020} and \cite{Cruz2020b} propose a coverage metric to quantify the sensing quality of bus fleet in different application scenarios, and found that: (1) The bus-based sensing has distinctive advantage over fixed-location (stationary) sensing for some urban applications; (2) The monitoring requirements have a significant impact on the sensing quality.
        
        In another line of research, the optimization of sensor deployment for bus-based urban sensing mainly focuses on assigning sensors to bus lines (depending on their spatial characteristics) or individual vehicles (based on their GPS trajectories). All these models can be characterized as subset selection problems, by selecting a subset of lines or vehicles to install sensors. \cite{YLL2012} adopt a chemical reaction optimization approach to select a subset of bus lines to install sensors. \cite{AD2017} design a greedy heuristic to allocate sensors to bus lines to maximize the total road length covered by the sensors. \cite{Kaivonen2020} employ route coverage image analyses to evaluate the percentage of area  covered by various bus route combinations. \cite{Saukh2012} design an evolutionary algorithm to select a subset of tram lines based on the service timetables to achieve satisfactory coverage. \cite{Gao2016} and \cite{WLQWL2018} design greedy heuristic algorithms to select a subset of buses to maximize the sensing quality based on their historical GPS trajectories. \cite{TRMSO2020} design a heuristic algorithm to allocate sensors to buses considering the importance level of different geographical regions. \cite{Agarwal2020} design an approximate algorithm to select a subset of hosts from a fleet of mixed buses and taxis based on their trajectories to maximize the spatial-temporal coverage.

        These studies, which are further summarized in Table \ref{tabSLP}, all employ heuristic, meta-heuristic, or approximate algorithms due to the size and  complexity of the problem. A more comprehensive review of sensor deployment issues in drive-by sensing can be found in \cite{JHL2023}.

        \begin{table}[h!]
		\centering
		\caption{Summary and comparisons of studies on the sensor deployment problem based on bus fleets.}
		\small{
		\begin{tabular}{|m{0.19\textwidth} | m{0.11\textwidth}  | m{0.24\textwidth} |m{0.1\textwidth}| m{0.2\textwidth}|}
		\hline
		Study & Monitoring scenario & Problem scale & Sensor host & Solution method
		\\
		\hline
		\cite{YLL2012} & Air quality & 91 bus lines in  Hong Kong Island, China & Bus lines           & Chemical reaction (meta-heuristic) 
		\\
		\hline
		\cite{AD2017} & Road potholes           & 713 bus lines in London, UK & Bus lines     & Greedy heuristic 
            \\
		\hline
		\cite{Kaivonen2020} & Air quality   & 21 bus lines in Uppsala, Sweden & Bus lines              & Image analysis algorithm
            \\
            \hline
            \cite{Saukh2012} & Air quality   & 13 tram lines in Zurich, Switzerland & Tram lines              & Evolutionary algorithm
            \\
		\hline
		\cite{Gao2016} & Air quality                       & Trajectories of 1,415 buses in Hangzhou, China & Vehicles   & Greedy heuristic    
		  \\
		\hline
		\makecell[l]{\cite{WLQWL2018}}    & --                   & Bus trajectories in Beijing, China            & Vehicles                                                          & Greedy approximate algorithm  
            \\
		\hline
		\cite{TRMSO2020}    & Heat island   & Trajectories of 20 buses in Athens, Georgia   & Vehicles   & Greedy heuristic algorithm  
            \\
            \hline
		\cite{Agarwal2020}    & Air quality   & Taxi and bus trajectories in San Francisco, USA   & Vehicles   & Approximation algorithm  
            \\
		\hline
		\makecell[l]{This study} & General & 167 bus lines in Chengdu, China, with 27,144 trips              & Trip chains & Linear integer programs
            \\
		\hline
		\end{tabular}
		}
		\label{tabSLP}
        \end{table}

                Different from these approaches, this paper introduces a trip-based solution framework for the mobile sensor deployment problem. By decoding the spatial-temporal information of individual service trips, we simultaneously determine trip chains and sensor assignment, which allows the sensor to be optimally circulated in the space-time domain with minimum interference to the level of bus service.

        \subsection{Minimum fleet in bus scheduling}
        At the core of our proposed methodology is a bus scheduling problem with minimum fleet size consideration. The minimization of fleet size is an important task in bus scheduling, and has been addressed in the literature using analytic modeling approaches \citep{Newell1971, DO2019}, graphical modeling approaches \citep{Ceder2016}, network flow approaches \citep{KMS2006, BK2009, LLX2019, Gkiotsalitis2023}, heuristic approaches \citep{Bartlett1957}, and bipartite graph approaches \citep{Levin1971}. It is generally agreed that the minimum fleet size problem, due to its combinatoric nature, has multiple optimal solutions. This means that multiple feasible configurations of vehicle trip chains exist that satisfy the minimum fleet size. A few studies have considered other objectives with the minimum fleet size as constraints. \cite{LC2020} develop integer programming models to minimize the number of required battery chargers for electric bus fleets while ensuring minimum fleet size. \cite{LJGC2022} optimize transfer coordination among different transit lines while complying with the minimum fleet size constraint. 
        
        In this work, the minimum fleet size issue is addressed using an efficient bipartite graph matching approach, and is integrated with the bus drive-by sensing framework as an optimization objective or a side constraint. It is demonstrated that exploring the non-uniqueness of the minimum fleet problem is highly effective in boosting the sensing power of bus fleets.

	\section{Problem statement and preliminaries}\label{secPS}
        
        \subsection{Problem description}
        
        \begin{definition}[\textbf{Trip-based mobile sensor deployment (MSD) problem for buses}]
            Consider a target monitoring area with a set of bus lines $\mathcal{L}$. Each line $l\in \mathcal{L}$ has a set $I_{l}$ of timetabled service trips. Each trip $i\in I_l$ is represented as a $4$-tuple ($p^i$, $q^{i}$, $t_{li}$, $\tau_i$), denoting departure terminal, arrival terminal, start time, and trip service time, respectively. The MSD problem considered in this paper aims to allocate a given number of sensors to the buses (such buses are called instrumented), and determine the set of trips to be carried out by each individual bus such that all timetabled trips in the system are fulfilled, while the quality of sensing provided by instrumented buses is maximized.
            \label{defMSD}
        \end{definition}
        
By definition, the MSD problem for buses needs to simultaneously determine the sequence of trips (coined {\it trip chain} in this paper) to be executed by a single bus, as well as the allocation of sensors to the buses. Due to the one-to-one correspondence of trip chain and bus, it suffices to associate sensors with trip chains, which will be the main approach taken in this paper. The following describes the context and assumptions of our problem. 

\begin{enumerate}
\item The target monitoring area is spatially meshed into grids while the time is discretized into intervals. The sensing coverage of a grid at certain time interval is calculated based on its intersection with the space-time trajectories of instrumented buses (buses equipped with sensors). Calculation details are presented in Section \ref{sec_sensing_reward}. 

\item  All bus lines in the target area have fixed timetables to be strictly followed. A given bus can serve multiple timetabled trips in sequence, subject to trip connection time constraints. For a given trip $i$ the service time $\tau_i$ is known, and the bus moves along the route at a constant speed (which could vary among different trips). The last assumption is for inferring the coverage of a grid during certain window. 

\item The size of the bus fleet to carry out the scheduled service trips and sensing tasks is to be determined within the model, following a minimum fleet size principle. This reflects our consideration of bus operations with limited interference from drive-by sensing tasks. 
\end{enumerate}

Despite these assumptions, the methodological framework and insights are generalizable. In fact, a few model extensions to accommodate more practical  considerations are discussed in Section \ref{subsecdiscussion}.

        \subsection{Quantifying the sensing quality}\label{sec_sensing_reward}
        
        The urban sensing requirements vary by applications. For example, air quality and traffic state monitoring require continuous and widespread sensing; noise pollution and heat island require less frequent but reliable sensing. A proper way to assess the sensing quality needs to consider three aspects: (1) spatial extent; (2) temporal frequency; and (3) temporal duration. Meanwhile, the quality of sensing can be interpreted as the number of data points collected by one or more sensing vehicles in an area during certain interval \citep{JHL2023}.

        To reflect these considerations, we mesh the target area into spatial grids $g\in G$, and discretize the time horizon into sensing intervals $t\in T$. The length of each sensing interval is denoted $\Delta$. Let $N_{gt}\in \mathbb{Z}_+$ be the number of sensors (sensing vehicles) within grid $g$ during $t$, we define the sensing reward for a given pair $(g, t)$ as
        \begin{equation}\label{eqndefngt}
        n_{gt}\doteq 
        \begin{cases}
        1,\qquad\hbox{if}~ N_{gt}\geq 1
        \\
        0,\qquad \hbox{if}~N_{gt}=0
        \end{cases}
        \end{equation} 
        Such a reward function is frequently used in the drive-by sensing literature \citep{AD2017, TRMSO2020}. In addition, different monitoring areas or times of the day may be assigned different sensing priorities, which are translated into spatial $w_g$ and temporal $\mu_t$ weights, such that $\sum_{g\in G}w_g=1$, $\sum_{t\in T}\mu_t=1$. Then, the overall sensing objective is given below:
        \begin{equation}\label{eqnPhidef}
        \max \Phi =\sum_{g \in G}w_{g}\sum_{t \in T}\mu_{t}n_{gt} 
        \end{equation}
\noindent In other words, $\Phi$ is a space-time weighted average of the rewards $n_{gt}$. Such a form addresses the abovementioned considerations as follows:
\begin{itemize}
\item[(1)] The extent of spatial and temporal coverage is directly reflected in the objective;

\item[(2)] Different sensing frequency requirements can be accommodated by setting the value of $\Delta$; a large $\Delta$ means a low sensing frequency requirement. 

\item[(3)] Using $n_{gt}$ instead of $N_{gt}$ in the objective prevents over-concentration of sensing vehicles in one place or one time, as the marginal gain is zero per \eqref{eqndefngt}.
\end{itemize}
        
        Intuitively, assuming equal weights $w_g$'s and $\mu_t$'s, the quantity $\Phi\in(0, 1)$ can be interpreted as the proportion of space-time pairs $(g,t)$ that are covered by at least one sensor. The objective \eqref{eqnPhidef} will be used in the remainder of the paper.

    \section{Model development for the mobile sensor deployment problem}\label{sec_MLP_two_mode}

We are given a set of bus lines $\mathcal{L}$ in the target area, and for each $l\in\mathcal{L}$, its timetabled trips are known and must be fulfilled. A bus in line $l$ can serve multiple trips in sequence, and the set of such trips is called a trip chain. 

\begin{figure}[h!]
\centering
\includegraphics[width=\textwidth]{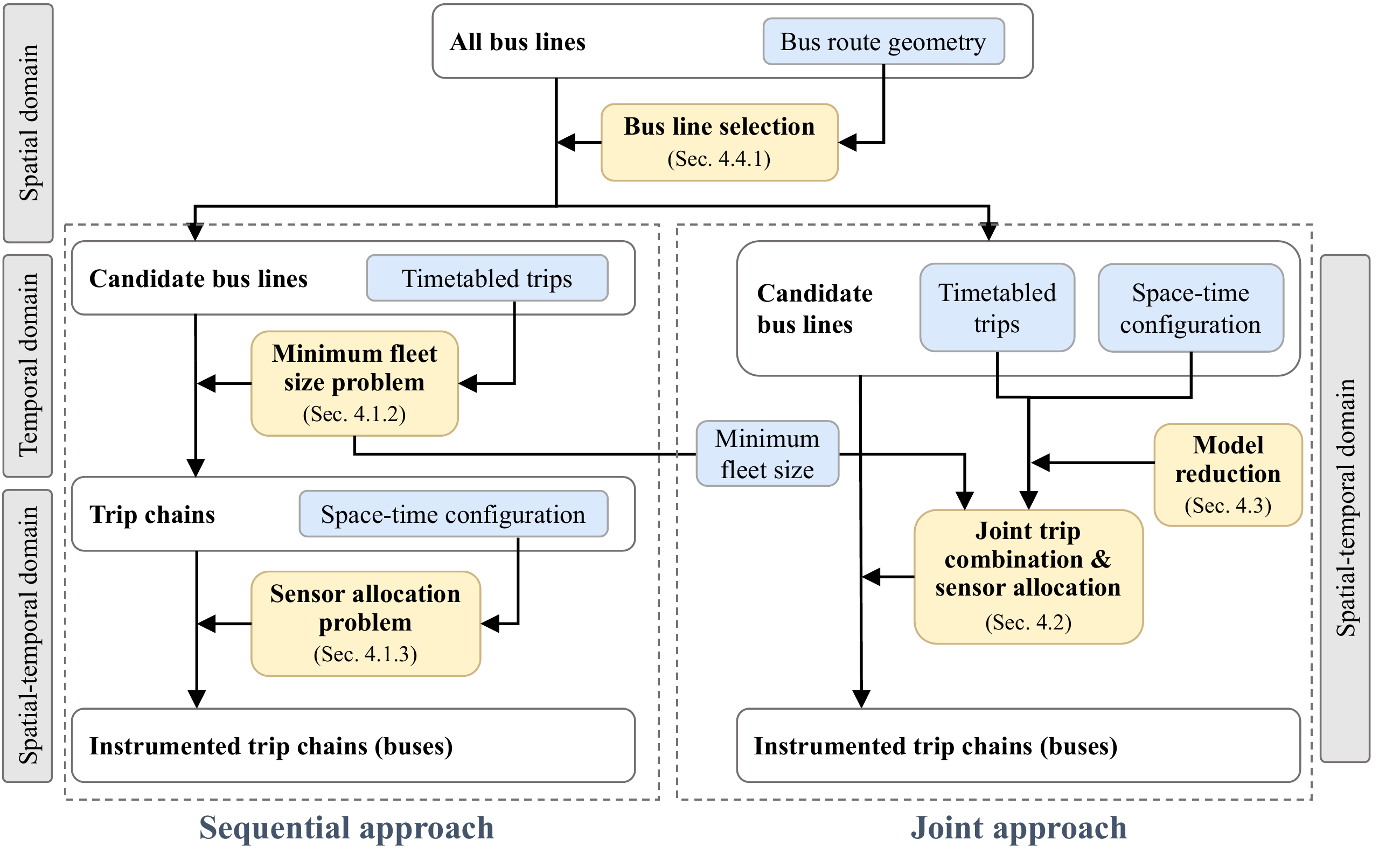}
\caption{Methodological framework for the multi-stage bus-based MSD problem.}
\label{figflow}
\end{figure}

This paper proposes a sequential and a joint solution approach for the bus-based mobile sensor deployment (MSD) problem, which are illustrated in Figure \ref{figflow}. Both approaches are multi-stage, relying on a pre-processing step that selects candidate bus lines for sensor allocation. The line pre-selection step is necessary as it not only reduces the candidate lines by a substantial amount (see Remark \ref{lsrmk}), but also effectively decouples the spatial and temporal aspects of the sensing objective, allowing subsequent modules to be solved within each bus line in a distributed manner (see Remark \ref{rmkbi}). 

Once the candidate lines are determined, the following procedures are executed:
\begin{itemize}
\item {\bf Sequential approach}: The sequential approach attempts to combine timetabled trips within a given line into trip chains via a bipartite graph matching approach \citep{Levin1971}, such that all the trips can be executed by a fleet of minimum size. Then,  sensor allocation problem is optimized as a subset selection problem, which determines whether a trip chain (the serving bus) is instrumented.

\item {\bf Joint approach}: The joint approach simultaneously combines trips  while assigning sensors to the resulting chains. For bus operational considerations, the minimum fleet size, which is obtained from the sequential approach, is imposed as a constraint. To further enhance the computational performance, a model reduction procedure is developed by analyzing the minimum fleet problem, which substantially reduces the number of binary decision variables and removes inefficient constraints. 
\end{itemize}

The joint approach explores the fact that the minimum fleet size problem has multiple trip chain solutions. Thus, simultaneously optimizing trips combination and sensor allocation allows us to find the minimum-fleet solution that promises the best sensing outcome.

Table \ref{tab_notations} lists key parameters and variables used in our model.

\setlength\LTleft{0pt}
	\setlength\LTright{0pt}
	\begin{longtable}{@{\extracolsep{\fill}}rl}
		\caption{Notations and symbols}
            \label{tab_notations} 
            \\
	    \hline
			\multicolumn{2}{l}{Sets}     \\
			\hline
$\mathcal{L}$ & Set of all bus lines in the target area;
\\
$L$    & Set of selected bus lines in the target area;
			\\
$T$   & Set of discrete time periods;
                \\
$G$  & Set of spatial grids reached by the bus line set $\mathcal{L}$ in the target area;
			\\
$I_{l}$ & Set of timetabled service trips of bus line $l\in L$;
\\
$U_l$ & Set of pull-out trips from the depot of line $l\in L$;
\\
$V_l$ & Set of pull-in trips to the depot of line $l\in L$;
                \\
$C$ & Set of trip chains;
                \\ 
$C_l$ & Set of trip chains in line $l\in L$;
                \\
$B_{l}$ & Set of buses in line $l\in L$.
                \\
                \hline
			\multicolumn{2}{l}{Parameters and constants}                                                                                                         			\\
			\hline
		$\delta_{gl}$ & Binary constant that equals $1$ if line $l$ intersects grid $g$, and $0$ otherwise;
		\\
		$N_G$ & Total number of grids covered by all lines in $\mathcal{L}$;
		\\
		$\gamma$ & Proportion of grids to be covered by the selected bus lines
		\\
		$N_{B_l}$ & Total number of buses in line $l\in L$;
		\\
                $N_{I_l}$    & Total number of timetabled trips of bus line $l \in L$ within the time horizon;
                \\
                $N_{S}$          & Total number of sensors;
                \\
                $\tau_i$           & End-to-end bus running time to serve a complete trip $i$;
			\\
                $t_{lij}$            & Bus deadheading time of line $l$, from the end of trip $i$ to the start of trip $j$;
                \\
                $t_{li}$            & Vehicle departure time of trip $i$ of bus line $l$
                \\
                $n_{igt}$           &  Binary constant that equals $1$ if trip $i$ covers grid $g$ in time period $t$;
                \\
                $w_{g}$  & Spatial sensing weight of grid $g$
                \\
                $\mu_{t}$  & Temporal sensing weight of time period $t$
                \\
                $q_{gt}^{m_l}$ & Binary parameter that equals $1$ if line $l$ covers grid $g$ during $t$, when the number of \\
                & sensors assigned to line $l$ is $m_l$.
                \\
			\hline
			\multicolumn{2}{l}{Auxiliary variables}                                                                                                \\
			\hline
             $x_g$ & Binary variable that equals $1$ if grid $g$ is covered by the selected lines;
             \\
                $N_{l}^{\text{min}}$         & The minimum bus fleet size of bus line $l$, $l \in L$
                \\
                              $\xi_{ci}$          & Binary variable that equals $1$ if trip chain $c$ includes trip $i$;
                \\
                $n_{gt}$  & Binary variable that equals $1$ if grid $g$ is covered by a sensor during $t$;
                \\
			\hline
			\multicolumn{2}{l}{Decision variables}                                                                                                 \\
			\hline
		$x_l$  & Binary variable that equals $1$ if line $l\in\mathcal{L}$ is selected for subsequent optimization;
		\\	
                $y_{lij}$  & Binary variable that equals $1$ if trip $j\in I_l$ is served after $i\in I_l$ with the same bus; 
                \\
                $z_{c}$  & Binary variable that equals $1$ if trip chain $c$ (the executing bus) is assigned a sensor;
                 \\
                $\xi_{cij}$         & Binary variable that equals $1$ if trip chain $c$ includes trip sequence $i$-$j$;
                \\
              $m_l$        & Number of sensors allocated to bus line $l\in L$;
                \\
                \hline
	\end{longtable}

    \subsection{The MSD problem with the sequential optimization approach}\label{sec_OS}
  We develop a three-stage optimization procedure: (1) bus line selection (Section \ref{sec_BR_selection}), (2) minimum fleet size determination (Section \ref{sec_MFS}), and (3) sensor allocation (Section \ref{sec_ASTV}).

 \subsubsection{Bus line selection}\label{sec_BR_selection}

         Given the set of all bus lines $\mathcal{L}$ in the target area, we aim to select a subset of lines $L\subset \mathcal{L}$ such that at least a proportion $\gamma$ of the spatial grids in the target area are covered by $L$.

        Let the binary parameter $\delta_{gl}=1$ if line $l\in\mathcal{L}$ intersects grid $g$, the binary decision variable $x_l=1$ if line $l$ is selected, and the binary auxiliary variable $x_g=1$ if grid $g$ is covered by at least one of the selected lines. The total number of spatial grids that can be covered by the bus line set $\mathcal{L}$ is $\sum_{g\in G}x_g$, which should be no less than $\gamma N_G$.  The mathematical model of bus line selection is as follows:
        \begin{equation}\label{FBRS1}
            \min_{x_l} \sum_{l\in \mathcal{L}}x_{l}
        \end{equation}
      \noindent such that
        \begin{eqnarray}
        \label{BLS1}
    		-M(1-x_{g})  \leq \sum_{l\in \mathcal{L}}\delta_{gl}x_l - 1  &&  \forall g \in G
    		\\
    		\label{FBRS3}
    		\sum_{l\in \mathcal{L}}\delta_{gl}x_l  \leq Mx_{g} &&  \forall g \in G
  \\
  \label{BLS2}
 \gamma N_G \leq \sum_{g\in G}x_g  &&
  \\
  \label{BLS3}
  x_l \in\{0, 1\} && \forall l\in \mathcal{L}
  \\
  \label{BLS4}
  x_g  \in\{0, 1\} && \forall g\in G
  \end{eqnarray}
  
 \noindent The objective function \eqref{FBRS1}  minimizes the total number of selected bus lines. Eqns \eqref{BLS1}-\eqref{FBRS3} express the binary state of coverage for grid $g$, and are equivalent to the following via the big-M method:
 $$
   x_g=\min\left\{1,\, \sum_{l\in \mathcal{L}}\delta_{gl}x_l\right\}\qquad\forall g\in G
 $$
 \noindent where $M$ is a sufficiently large number. \eqref{BLS2} ensures that the selected bus lines cover $N_G$ grids. In summary, the bus line selection sub-problem is formulated as a linear integer program \eqref{FBRS1}-\eqref{BLS4}.
 
 \begin{remark}
 In case $\gamma=1$, the bus line selection problem can be simplified as a set covering problem:
 \begin{eqnarray}
 \label{blpseqn1}
 \min \sum_{l\in\mathcal{L}}x_l
 \\
  \label{blpseqn2}
 \sum_{l\in\mathcal{L}}\delta_{gl}x_l\geq 1 && \forall g\in G
 \\
  \label{blpseqn3}
 x_l\in \{0,\,1\} && \forall l\in\mathcal{L}
 \end{eqnarray}
 \noindent In our numerical case studies, we set $\gamma=1$ and use the simpler set covering formulation.
 \end{remark}

 The purpose of the line selection sub-module is two fold: (1) to limit the number of bus lines considered in subsequent optimization procedures to reduce the computational complexity; (2) to spatially decouple the sensing tasks by individual lines such that the MSD can be solved in a distributed manner (also see Remark \ref{rmkbi}).

    \subsubsection{Minimum fleet size problem}\label{sec_MFS}

    In this stage, we aim to determine the minimum bus fleet size required to fulfill all timetabled trips within each $l\in L$. This is equivalent to a bipartite graph maximal matching problem \citep{Levin1971}. As a result, we obtain a partition of trips, with each subset called a trip chain.  

    \begin{equation}\label{FMFS1}
        \min_{y_{lij}} \sum_{l \in L}N_l^{\text{min}}=\sum_{l\in L} \left( N_{I_l} - \sum_{i,j\in I_l} y_{lij} \right)
    \end{equation}

    \begin{eqnarray}
		\label{FMFS2}
		t_{lj}-(t_{li}+\tau_{i})-t_{lij} \geq -M(1-y_{lij}) & & \forall l \in L , \forall i,\,j \in I_{l}
		\\
		\label{FMFS3}
		\sum_{i\in I_{l}}y_{lij} \leq 1 & & \forall l \in L , \forall j \in I_{l}
		\\
		\label{FMFS4}
		\sum_{j \in I_{l}}y_{lij} \leq 1 & & \forall l \in L , \forall i \in I_{l}
		\\
		\label{FMFS6}
		y_{lij} \in\{0,\,1\} & & \forall l \in L , \forall i,\,j \in I_{l}
    \end{eqnarray}

\noindent The objective \eqref{FMFS1} is to minimize the required bus fleet size for all selected bus lines $l \in L$. Constraint \eqref{FMFS2} indicates whether trip $j$ of line $l$ can be conducted by the same vehicle after serving trip $i$, which is equivalent to the following:
    \begin{equation}\label{eqnconn}
    y_{lij}=
    \begin{cases}
    0~\hbox{or}~1\quad &\text{if}~~t_{li}+\tau_i+t_{lij}\leq t_{lj}
    \\
    0 \quad &\text{if}~~t_{li}+\tau_i+t_{lij}> t_{lj}
    \end{cases}
    \end{equation}

 \noindent Constraints \eqref{FMFS3} and \eqref{FMFS4} ensure that each trip can be connected with no more than one predecessor and successor trips, respectively.

The idea of the minimum fleet size problem is straightforward: By combining trips in a tail-to-head fashion, subject to the connectivity constraints \eqref{FMFS3}-\eqref{FMFS4}, we form a trip chain, which is executed by a single bus. The minimum fleet size is thus obtained by maximizing such connections (i.e. $y_{lij}$'s). Upon solving the problem \eqref{FMFS1}-\eqref{FMFS6}, the resulting trip chains form the set $C$, and every $c\in C$ is an ordered set of trips $c=\{i_1, i_2,\ldots \}$.

    \subsubsection{Allocating sensors to trip chains} \label{sec_ASTV}

   In the final stage, the model allocates a given number of sensors to the  trip chains generated in the minimum fleet size problem, as each trip chain is accommodated by the same vehicle. The goal is to maximize the sensing quality quantified as the reward:

    \begin{equation}\label{FOS1}
        \max_{z_c} \sum_{g \in G}w_{g}\sum_{t \in T}\mu_{t}n_{gt} 
    \end{equation}

    \begin{eqnarray}
		\label{FOS2}
		n_{gt} \leq \sum_{c \in C}\sum_{i \in c}z_{c} n_{igt} & & \forall g \in G, \forall t \in T
		\\
		\label{FOS3}
		\sum_{c \in C}z_{c} \leq N_{S}
		\\
		\label{FOS4}
		z_{c} \in \{0,\,1\} & & \forall c \in C
            \\
            \label{FOS5}
            n_{gt} \in \{0,\,1\} & & \forall g \in G, \forall t \in T
    \end{eqnarray}

\noindent The objective function \eqref{FOS1} maximizes the sensing reward; constraint \eqref{FOS2} indicates whether grid $g$ is covered by the sensing vehicles in time period $t$; constraint \eqref{FOS3} ensures that the number of installed sensors does not exceed the total. Note that in \eqref{FOS2}, $n_{igt}$'s are binary constants. 

\begin{remark}
Constraint \eqref{FOS2} should be more accurately written as
$$
n_{gt}=\min\left\{1,\, \sum_{c \in C}\sum_{i \in c}z_{c} n_{igt}\right\}
\qquad \forall g\in G,\,\forall t\in T
$$
In fact, one can easily show that this is equivalent to \eqref{FOS2} given that $n_{gt}\in \{0,1\}$ and the objective is to maximize a convex combination of $n_{gt}$'s. 
\end{remark}

\begin{remark}\label{lsrmk}
The bus line selection module, as a preprocessing step of the proposed sequential procedure, can considerably reduce the decision space. Taking the Chengdu bus network as an example (see Section \ref{subsecmodel} for more details), where there are 167 bus lines in the study area. Without line pre-selection \eqref{blpseqn1}-\eqref{blpseqn3}, the minimum fleet size is 2824, meaning that there are 2824 trip chains to be considered for sensor deployment, while this number is reduced to 684 with line pre-selection (38 lines in total). The subsequent impact on computational efficiency and solution optimality is analyzed in Section \ref{subsubsecEffls}.
    
    \end{remark}

    \subsection{The MSD problem with the joint optimization approach}\label{sec_AS}
    In the previous section, the formation of trip chains and their sensor allocation are determined sequentially. Due to the non-uniqueness of the solution to the minimum fleet size problem, $N_l^{\text{min}}$ could correspond to multiple trip chain configurations, which render different sensing outcomes in the subsequent sensor allocation problem. As there are no easy ways to systematically enumerate all such possible trip chain configurations, we develop a joint approach that aims to maximize the sensing quality while treating the minimum fleet size as a side constraint. This tactic allows us to explore the coordination between trip chain formation (i.e. bus scheduling) and sensing tasks. 
    
\begin{remark}\label{rmkbi}
    It should be noted that the minimum fleet size sub-module can be performed independently for each line, as both the objective and the constraints of \eqref{FMFS1}-\eqref{FMFS6} are entirely decoupled by $l$. On the other hand, the sensing reward offered by different lines are inter-dependent, since a timed grid $(g, t)$ may be covered by sensors from different lines that overlap over $g$. Fortunately, in the line selection stage, such overlap is restricted by minimizing the number of lines to cover as many grids as possible. Therefore, it makes sense to decouple the joint optimization approach by individual bus lines to reduce the computational complexity. 
    \end{remark}
    
    Building on Remark \ref{rmkbi}, we propose a bi-level approach:
    
 \begin{eqnarray}
 \label{eqnupper}
 \text{(Upper Level)} & & \max_{\left\{m_l\, : \, l\in L\right\}} \Phi
 \\
 \label{eqnlower}
 \text{(Lower Level)} & & \max_{C_l,\, \left\{z_c\, : \, c\in C_l\right\} \,\vert m_l} \Phi_l \qquad \forall l\in L
 \end{eqnarray}
\noindent where \eqref{eqnupper} determines the number of sensors $m_l$ allocated to each line $l\in L$, and \eqref{eqnlower} jointly solves for trip chain configuration $C_l$ as well as sensor assignment $\{z_c: c\in C_l\}$, assuming  $m_l$ is given. The upper-level objective is the sensing reward over the entire target area, while the lower-level sensing objective is restricted to the coverage area of line $l$. Such a decomposition scheme allows us to solve the joint trip-chain-configuration-and-sensor-allocation problem within each line $l$ independently.

   \subsubsection{The lower-level problem} \label{sec_AS_SL}
   
	The lower-level problem simultaneously determines the trip chain set $C_l$ and sensor allocation $z_c: c\in C_l$ for each bus line $l\in L$,  subject to the minimum fleet size constraint. The minimum fleet size $N_l^{\text{min}}$ is calculated according to Section \ref{sec_MFS}.

	For reason that will become clear in Remark \ref{rmkreason}, the joint approach needs to consider non-service trips of being pulled-out from the depot (denoted $u\in U_l$) and pulled-in to the depot (denoted $v\in V_l$). Without loss of generality, the trips $u$'s and $v$'s can be treated as normal service trips when we express trip connection and flow conservation constraints. Note that the trips from/to the depot do not contribute to the sensing objective, as we do not wish to dive into details of the depot locations and pull-in/out trips. The mathematical formulation of the lower-level problem  is as follows: $\forall l \in L$,

   \begin{equation}\label{FAS1}
        \max \Phi_l= \sum_{g \in G}w_{g}\sum_{t \in T}\mu_{t}n_{gt}
    \end{equation}
    \begin{eqnarray}
        \label{FAS2}
	   t_{lj}-(t_{li}+\tau_{i})-t_{lij} \geq -M(1-\xi_{cij}) & &  \forall i,\,j \in I_{l}, \forall c\in C_l
	\\
        \label{FAS3}
        \sum_{c\in C_l}\sum_{i \in U_l\cup  I_{l}}\xi_{cij} = 1 & & \forall j \in I_{l}
        \\
        \label{FAS35}
     \sum_{c\in C_l} \sum_{j \in I_{l}\cup V_l}\xi_{cij} =1   & & \forall i \in I_{l}
        \\
        \label{FAS4}
        \sum_{u\in U_l}\sum_{j \in I_{l}}\xi_{cuj} = 1 & & \forall c \in C_{l}
        \\
        \label{FAS5}
        \sum_{i \in I_{l}}\sum_{v\in V_l} \xi_{civ} = 1  & & \forall c \in C_{l}
        \\
        \label{FAS6}
        \sum_{i \in I_{l}\cup U_l }\xi_{cij} = \sum_{i \in  I_{l}  \cup V_l }\xi_{cji} & & \forall c \in C_{l}, \forall j \in I_{l}
                \\
        \label{FAS8}
        \xi_{ci} = \sum_{j \in I_{l}\cup V_l}\xi_{cij} & & \forall c \in C_{l}, \forall i \in  I_{l}
        \\
        \label{FAS7}
        N_{I_l}- \sum_{c\in C_l}\sum_{i,j \in I_{l}}\xi_{cij} = N_{l}^{\text{min}}
        \\
		\label{FAS9}
		n_{gt} \leq \sum_{c\in C_l}\sum_{i \in I_l}z_{c} \xi_{ci} n_{igt} & & \forall g \in G, \forall t \in T
		\\
		\label{FAS10}
		\sum_{c \in C_{l}}z_{c} \leq m_l
		\\
		\label{FAS11}
		\xi_{cij} \in \{0,\,1\} & & \forall c \in C_{l}, \forall i,j \in I_{l}\cup U_l \cup V_l
            \\
            \label{FAS12}
		\xi_{ci} \in \{0,\,1\} & & \forall c \in C_{l}, \forall i \in I_{l}
            \\
		\label{FAS14}
		z_{c} \in \{0,\,1\} & & \forall c \in C_{l}
            \\
            \label{FAS15}
        n_{gt} \in \{0,\,1\} & & \forall g \in G, \forall t \in T
    \end{eqnarray} 

\noindent  The objective function \eqref{FAS1} maximizes the weighted sensing reward. Constraint \eqref{FAS2} expresses the trip connectivity constraint that is specific to trip chain $c$. Constraints \eqref{FAS3}  (\eqref{FAS35}) ensures that each trip belongs to exactly one trip chain and has exactly one predecessor (successor). Constraints \eqref{FAS4} (\eqref{FAS5}) ensures that each trip chain $c$ begins (ends) with a pull-out (pull-in) trip. Constraint \eqref{FAS6} expresses flow balance upstream and downstream of a trip $j$. Constraint \eqref{FAS7} ensures that the fleet size under consideration for line $l$ is minimal, where $N_l^{\text{min}}$ is obtained from the minimum fleet size sub-problem \eqref{FMFS1}-\eqref{FMFS6}. \eqref{FAS8} is the trip consistency constraints, which states that a trip belongs to a trip chain $c$ if and only if it has one successor in $c$. Inequalities \eqref{FAS9} indicate whether spatial grid $g$ is covered during time $t$ by at least one sensor. Constraint \eqref{FAS10} limits the number of sensors to be no greater than the total allocated to bus line $l$, where $m_l$ is an exogenous parameter to this problem, which comes from the upper-level problem.

\begin{remark}\label{rmkreason}
The reason for introducing the dummy trips from/to the depot is apparent from \eqref{FAS3}-\eqref{FAS8}, where the trip chain balance and consistency principles would fail without $u$'s or $v$'s. Such notions are not necessary in the sequential approach because the minimum fleet problem only determines the connectivity between two trips (i.e. $y_{lij}\in \{0,1\}$), and the trip chains are only formed upon solving the optimization problem. On the other hand, to explicitly accommodate trip chains using equalities and inequalities within an optimization formulation requires a flow-based approach, which entails flow balance and consistency constraints.
\end{remark}

    \subsubsection{The upper-level problem}\label{subsubsecUpper}

    The upper- and lower-level problems interact with each other through the variable $m_l\in\mathbb{Z}_+$, which is the number of sensors allocated to line $l\in L$. Upon solving \eqref{FAS1}-\eqref{FAS15} with $m_l$ sensors, the lower level returns the binary constant $q_{gt}^{m_l}$, which is equal to $1$ if line $l$ covers grid $g$ in time $t$. Based on such interaction, the upper-level problem reads:

    \begin{equation}\label{FSTR1}
        \max_{m_l: l\in L} \Phi=\sum_{g \in G}w_{g}\sum_{t \in T}\mu_{t}n_{gt}
    \end{equation}

    \begin{eqnarray}
        \label{FSTR2}
        n_{gt} \leq \sum_{l \in L}q_{gt}^{m_l} & & \forall g \in G, \forall t \in T
        \\
        \label{FSTR3}
        \sum_{l \in L} m_l \leq N_{S}
        \\
        \label{FSTR4}
        n_{gt} \in \{0,1\} & & \forall g \in G, \forall t \in T
                \\
        \label{FSTR5}
        m_l \in \mathbb{Z}_+ & & \forall l \in L 
    \end{eqnarray}  

 \noindent The objective function \eqref{FSTR1} maximizes the sensing reward in the target area; constraint \eqref{FSTR2} indicates whether spatial grid $g$ is covered by the sensing vehicles in time period $t$; \eqref{FSTR3} ensures that the total number of deployed sensors cannot exceed $N_S$. 
 
 The feasible region of the decision variables $m_l$, given by \eqref{FSTR3} and \eqref{FSTR5}, is a simplex in a $|L|$-dimensional space. To reduce the size of the feasible region, we notice that for each line $l\in L$ the sensing reward \eqref{FAS1} increases with $m_l$, but stalls past certain threshold $m_l=K_l$. Therefore, it suffices to replace \eqref{FSTR5} with 
 \begin{equation}\label{FSTR5'}
 0\leq m_l\leq K_l\qquad \forall l\in L
 \end{equation} 
\noindent Given that each grid only needs to be covered once in a sensing interval, it is expected that the number of sensors to saturate a line $K_l$ is small, which is indeed the case as confirmed in Section \ref{subsubsectionsaturate}.

Given the above discussion, we derive the following algorithm to solve the upper-level problem.
 
 \begin{algorithm}[H]
    \caption{Solving the upper-level problem}
                \label{Aupper}
    \hspace*{0.02in} {\bf Input:}
    The total number of sensors $N_S$; 
     \\
    \hspace*{0.02in} 
    {\bf Output:} Optimal line assignment solution $m_l:\,l\in L$
    \\
    \hspace*{0.02in} 
    {\bf Step 1:} Solve for the upper bound $K_l$, $l\in L$   
    \begin{algorithmic}[1]
    \State $a_1 \leftarrow -1$, $a_2 \leftarrow 0$, $m_l \leftarrow 0$;
        \While{$a_1<a_2$} 
        \State $m_l=m_l+1$;
        \State  Solve the lower-level problem \eqref{FAS1}-\eqref{FAS15} with $m_l$ sensors, return $\Phi_l$ and $q_{gt}^{m_l}$;
         \State $a_1 \leftarrow a_2$, $a_2 \leftarrow\Phi_l$;
        \EndWhile
        \State\Return $K_l = m_l-1$, $\big\{ q_{gt}^{m_l}: m_l=1, \ldots, K_l\big\}$.
    \end{algorithmic}
     \hspace*{0.02in} 
    {\bf Step 2:} Solve the main program \eqref{FSTR1}-\eqref{FSTR4} and  \eqref{FSTR5'}.
\end{algorithm}

    \subsection{Model reduction procedure}\label{subsecmr}
    
    In this section, we develop a few preprocessing techniques to further reduce the complexity of the optimization problem, by reducing the feasibility region of the trip connectivity constraint. Similar techniques have been seen in the construction of the vehicle-shareability network \citep{Vazifeh2018}. 

{\bf Step 1. Feasible trip pair.} In the original model, two trips $i, j\in I_l$ can form a feasible pair (that is, $\xi_{cij}=1$ for some $c\in C_l$) if
\begin{equation}\label{eqnconnectible}
0\leq t_{lj}-\big(t_{li}+\tau_i+t_{lij} \big),
\end{equation}
which means the bus arrives at the starting terminal $p^j$ of trip $j$ before the starting time of $j$. The right hand side of \eqref{eqnconnectible} is called the terminal idle time. On the other hand, it is reasonable to expect that the idle time should not be too long, as the minimum fleet size solution promotes high vehicle utilization. Therefore, we may impose the following artificial constraint:
\begin{equation}\label{eqnarti}
t_{lj}-\big(t_{li}+\tau_i+t_{lij} \big)\leq \delta
\end{equation}
\noindent where $\delta>0$ is the maximum idle time at the terminal. Note that as $\delta$ decreases, more feasible trip pairs will be eliminated, and the minimum fleet size could increase. Therefore, to ensure that the solution is not altered by this artificial constraint, one should choose $\delta$ carefully to guarantee the same $N_l^{\text{min}}$ as before.

{\bf Step 2. Selection of $\delta$}. Appropriate values of $\delta$ can be obtained through trial and error. Since each trial entails solving the problem \eqref{FMFS1}-\eqref{FMFS6} and \eqref{eqnarti}, which may be time-consuming, to limit the number of trials we employ a binary search algorithm within the interval $[0,\,\delta_{0}]$. Here, $\delta_0$ is an initial guess such that the resulting minimum fleet is equal to $N_l^{\text{min}}$. As a rule of thumb, $\delta_0$ can be taken as the duration to complete a round trip. The reason is that any solution to the original minimum fleet size problem \eqref{FMFS1}-\eqref{FMFS6} is unlikely to have a feasible trip pair $i,j$ with terminal idle time longer than $\delta_0$. Even if it does, the idling bus can be used to execute a round trip and still be in time to serve $j$, which leads to a new solution with fleet size $\leq N_l^{\text{min}}$, and all idle times below $\delta_0$.

{\bf Step 3. Model reduction}. We start with the set of all trip pairs $I_l\times I_l$, and use \eqref{eqnconnectible} and \eqref{eqnarti} to eliminate infeasible ones, the resulting set is denoted $\Omega_{\text{FTP}}$. Then, the trip connectivity constraints \eqref{FAS2} can be removed from the model, as we have eliminated $y_{lij}$'s or $\xi_{cij}$'s that must be zero. 

Algorithm \ref{A_1} executes the aforementioned steps to eliminate binary variables in the linear integer programs. 
\begin{algorithm}[H]
    \caption{Reducing binary variables for line $l\in L$}
                \label{A_1}
    \hspace*{0.02in} {\bf Input:}
    Trip set $I_{l}$; initial guess $\delta_0$ and maximum iteration $K$ for the binary search; 
     \\
    \hspace*{0.02in} 
    {\bf Output:} The set of feasible trip pairs $\Omega_{\text{FTP}}$
    \\
    \hspace*{0.02in} 
    {\bf Step 1:} Binary search for $\delta$ that guarantees $N_{l}^{\text{min}}$   
    \begin{algorithmic}[1]
    \State $\eta_l \leftarrow 0$; $\eta_r \leftarrow \delta_0$;
        \For{$k=1: K$} 
        \State  $\delta_k={1\over 2}(\eta_l+\eta_r)$
        \State Solve the problem \eqref{FMFS1}-\eqref{FMFS6} and \eqref{eqnarti} with maximum idle time $\delta_k$
            \If {the resulting fleet size $>N_{l}^{\text{min}}$}
                \State $\eta_l=\delta_k$;
                \Else 
                \State $\eta_r=\delta_k$;
            \EndIf
        \EndFor
        \State\Return $\delta \in [\eta_r,\, \delta_0]$
    \end{algorithmic}
     \hspace*{0.02in} 
    {\bf Step 2:} Reduce binary variables 
    \begin{algorithmic}[1]
    \For{$\forall i, j\in I_l$}
    \If {$t_{lj}-\big(t_{li}+\tau_i+t_{lij}\big)\notin [0,\, \delta]$}
    \State $y_{lij}\leftarrow 0$; $\xi_{cij}\leftarrow 0$;
    \EndIf
    \EndFor
    \end{algorithmic}
\end{algorithm}

To demonstrate the effectiveness of the proposed scheme, Table \ref{tab_MRP_3} shows the number of binary variables $\xi_{cij}$'s in a few example bus lines. The values of $\delta=100, 200, 300$ are tested, which reduce the binary variables by at least 70\%. Moreover, the big-M constraint \eqref{FAS2}, although conveniently converting logical statements into inequalities without introducing additional binary variables, could increase the solution space of the linear relaxation sub-problem, which could cause slow convergence when employing the branch-and-bound procedure. Therefore, the proposed model reduction technique not only reduces the binary variables, but also avoids the big-M constraints, which promises significant acceleration of the computational procedure. The actual computational times will be presented in Table \ref{tabCPT} in the numerical case study, which also shows that the optimal objective values are not affected by the model reduction procedure.

\begin{table}[h!]
    \centering
    \caption{Results of variables reduction in a few bus lines. The percentages represent variable reduction compared to the original problem.}
    \resizebox{\textwidth}{!}{
        \begin{tabular}{ccccccc}
        \hline
            \multirow{3}{*}{\begin{tabular}[c]{@{}c@{}}Bus\\lines\end{tabular}} & \multirow{3}{*}{\begin{tabular}[c]{@{}c@{}} $N_{I_l}$ \end{tabular}} & \multirow{3}{*}{\begin{tabular}[c]{@{}c@{}}$N_{l}^{\text{min}}$ \end{tabular}} & \multicolumn{4}{c}{Number of variables $\xi_{cij}$}            \\
            \cline{4-7}
             & &    & \multirow{2}{*}{Original} & \multicolumn{3}{c}{Algorithm \ref{A_1}}  \\ 
             \cline{5-7}
             &    &    &    & $\delta=300$ & $\delta=200$ & $\delta=100$      \\
             \hline
            \#125    & 160    & 26   & $26 \times 250^2$    & $26 \times 6602$ (-89.44\%)  & $26 \times 4844$ (-92.25\%)  & $26 \times 2632$ (-95.79\%)                  \\
            \#85   & 160    & 22   & $22 \times 160^2$    & $22 \times 6466$ (-74.74\%)   & $22 \times 4704$ (-81.63\%) & $22 \times 2612$ (-89.80\%)                  \\
            \#10   & 198    & 20   & $20 \times 198^2$    & $20 \times 10702$ (-72.70\%)  &  $20 \times 7826$ (-80.04\%) & $20 \times 4304$ (-89.02\%)                \\
            \#128   & 160    & 16   & $16 \times 160^2$    & $16 \times 6526$ (-74.51\%)   & $16 \times 4710$ (-81.60\%) & $16 \times 2592$ (-89.88\%)  \\
            \hline
        \end{tabular}%
    }
\label{tab_MRP_3}
\end{table}

    \subsection{Discussion of model extensions}\label{subsecdiscussion}
The proposed trip-based MSD problem is based on a simplification of real-world bus operations, as our main goal is to bridge the gap between bus operations and drive-by sensing. In this part, we briefly discuss a few model extensions to accommodate more realistic issues.

\begin{itemize}
\item {\bf Service gaps.} The service of a bus may be temporarily suspended for refueling/charging or staff rest/change. Assuming these service gaps take place during certain time of the day, we introduce a set of dummy trips $D_l$ for line $l$, and each trip $i\in D_l$ is expressed as the 4-tuple (see Definition \ref{defMSD}) $(p^{\text{sg}}, p^{\text{sg}}, t_{li}, \tau_i)$ where $p^{\text{sg}}$ is the terminal at which the bus is parked during the service gap, $t_{li}$ is the start time, and $\tau_i$ is the duration. Then, the expression
$$
\sum_{i\in D_l}\xi_{ci}=1\qquad \forall c\in C_l
$$
\noindent enforces such a service gap for all trip chains. We note that different dummy sets can be created to accommodate different types of service gaps. 

\item {\bf Bus relocation.} So far in this paper, the connectivity of two trips $i$ and $j$ is governed solely by the temporal feasibility constraint \eqref{eqnconn}, which does not exclude the possibility that a bus needs to relocate to a different terminal after serving $i$ in order to carry out $j$. Such a relocation process obviously brings additional operational costs. To avoid bus relocation between any two trips $i, j\in I_l$, we simply set $t_{lij}$ to be a very large number, such that $y_{lij}=0$ per \eqref{eqnconn}. Obviously, limiting bus relocations would increase the minimum fleet size. In a more realistic setting where relocations are permitted and evaluated in a more comprehensive cost structure, they can be penalized in the minimum fleet size problem; see below. 

\item {\bf Bus operational costs.} In our paper, the overall optimization objective only concerns with sensing quality. In reality, bus operational costs include set-up costs, which are usually proportional to the fleet size, as well as running costs, which are comprised of service and relocation costs. These can be accommodated in the MSD problem by defining
$$
\Psi_1= \alpha N_l^{\text{min}} + \sum_{l\in L}\sum_{i,j\in I_l} \beta t_{lij} y_{lij} 
$$
\noindent in the sequential optimization, or
$$
\Psi_2= \alpha N_l^{\text{min}} + \sum_{l\in L}\sum_{c\in C_l}\sum_{i,j\in I_l}\beta t_{lij}\xi_{cij} 
$$
in the joint optimization. The first terms on the RHS denote the set-up costs, and the second terms represent relocation costs; $t_{lij}$ denotes the deadheading (relocation) time from the end of trip $i$ to the start of trip $j$, $\alpha$ and $\beta$ represent monetary costs per unit entity. $\Psi_1$ or $\Psi_2$ can be evaluated alongside the sensing objective $\Phi$ to form multi-objective optimization, which is beyond the scope of this work.

\item {\bf Uncertain service time.} We assume fixed service time $t_{li}$ for trip $i\in I_l$, and constant bus speed along the route. This is an ideal simplification, which could affect real-world sensing outcome when the bus is subject to random en-route delays. In general, this is not considered a significant issue as far as grid coverage is concerned, since the grid will be covered sooner or later, and a temporal misplacement could be remedied using data processing techniques such as numerical extrapolation. 

\end{itemize}

\section{Case study}\label{secCS}

    \subsection{Model setup and parameters}\label{subsecmodel}

\begin{figure}[h!]
    \centering
    \includegraphics[width=0.8\textwidth]{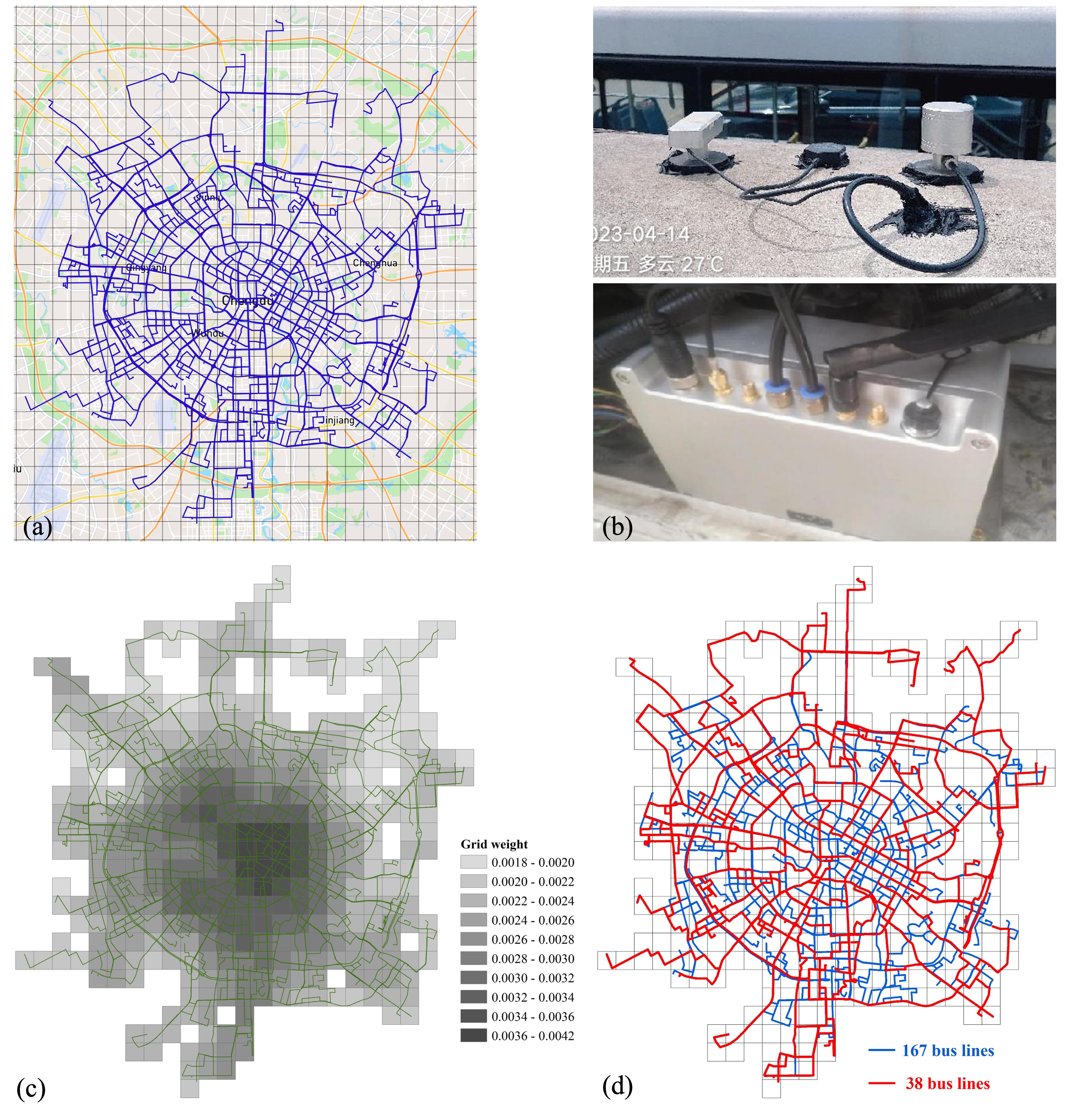}
    \caption{The target sensing area in Chengdu: (a) all 167 bus routes and spatial mesh (1km-by-1km); (b) mobile air quality sensor mounted on a bus;(c) spatial sensing weights $w_{g}$ of all the grids covered by the bus network; (d) 38 bus lines obtained from the bus line pre-selection procedure.}
    \label{fig_network}
    \end{figure}

    We consider an air quality monitoring scenario in Chengdu, China. The monitoring area is within the 4th Ring Road, with 400  spatial grids (1km$\times$1km) covered by 167 bus lines as shown in Figure \ref{fig_network}(a). Figure \ref{fig_network}(c) shows the grid weights $w_g$ in the sensing objective, which are calculated based on the distribution of main emission sources (including construction sites, factories and road traffic) and population density, which reflect the priority of environmental management. For $\gamma=1$ (all grids reachable by the 167 bus lines should be covered by the selected lines), the bus line pre-selection procedure \eqref{blpseqn1}-\eqref{blpseqn3} yields a total of 38 bus lines for subsequent optimization, which are shown in Figure \ref{fig_network}(d). 
    
    All the bus line geometries, including coordinates, terminal locations and directions, are obtained from the Amap Open Platform, where Amap.com is a digital map and navigation service provider. Bus trip service times are obtained from the Amap mobile App, and the bus deadheading times are estimated from shortest path search. The bus dispatch frequencies and timetables are obtained from a mobile App called Chelaile, which provides real-time bus positioning information.

\begin{figure}[H]
    \centering
    \includegraphics[width=\textwidth]{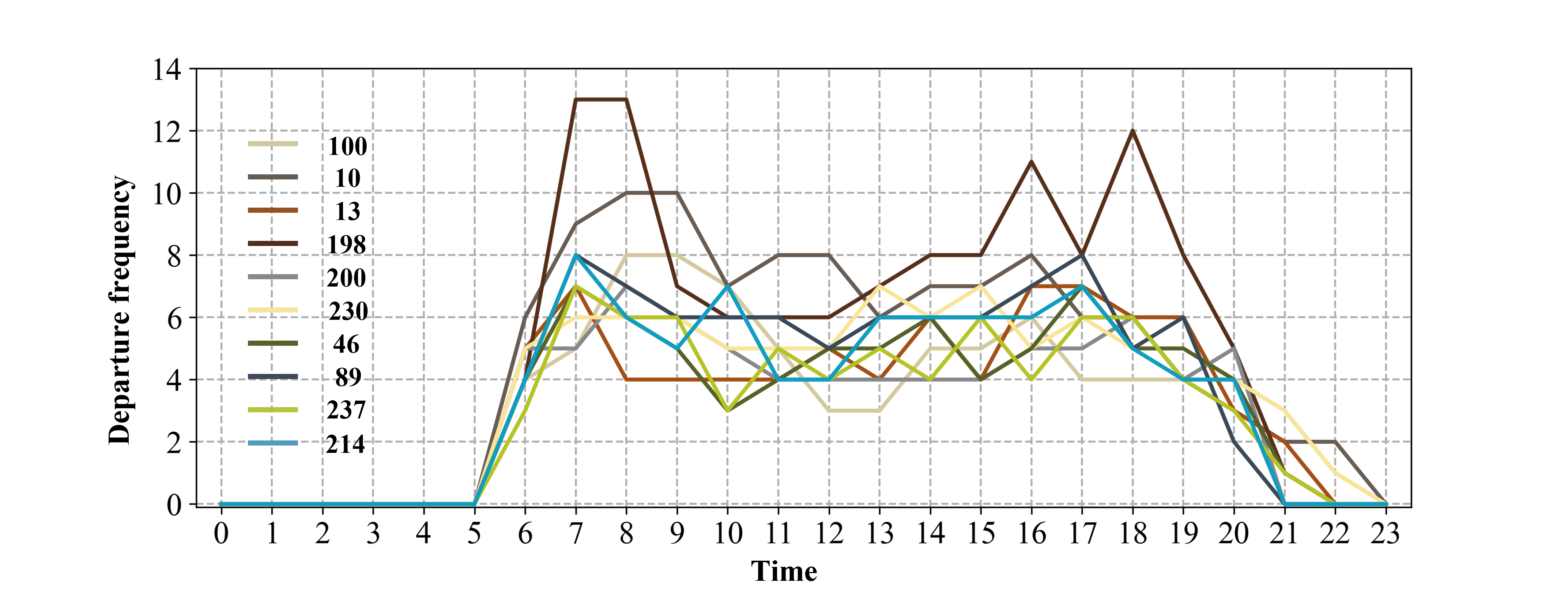}
    \caption{Departure frequency (veh/hour) of a few bus lines in the monitoring area}
    \label{fig_DF}
    \end{figure}

     The time horizon is set to be 7:00-22:00, spanning 15 hours. As shown in Fig. \ref{fig_DF}, this is consistent with the operating hours of most bus lines. Regarding the time interval for calculating sensing rewards, we consider three cases $\Delta=60, 90, 120$ (min). The corresponding temporal weights $\mu_t$'s are set to be equal.

\begin{figure}[H]
    \centering
    \includegraphics[width=1\textwidth]{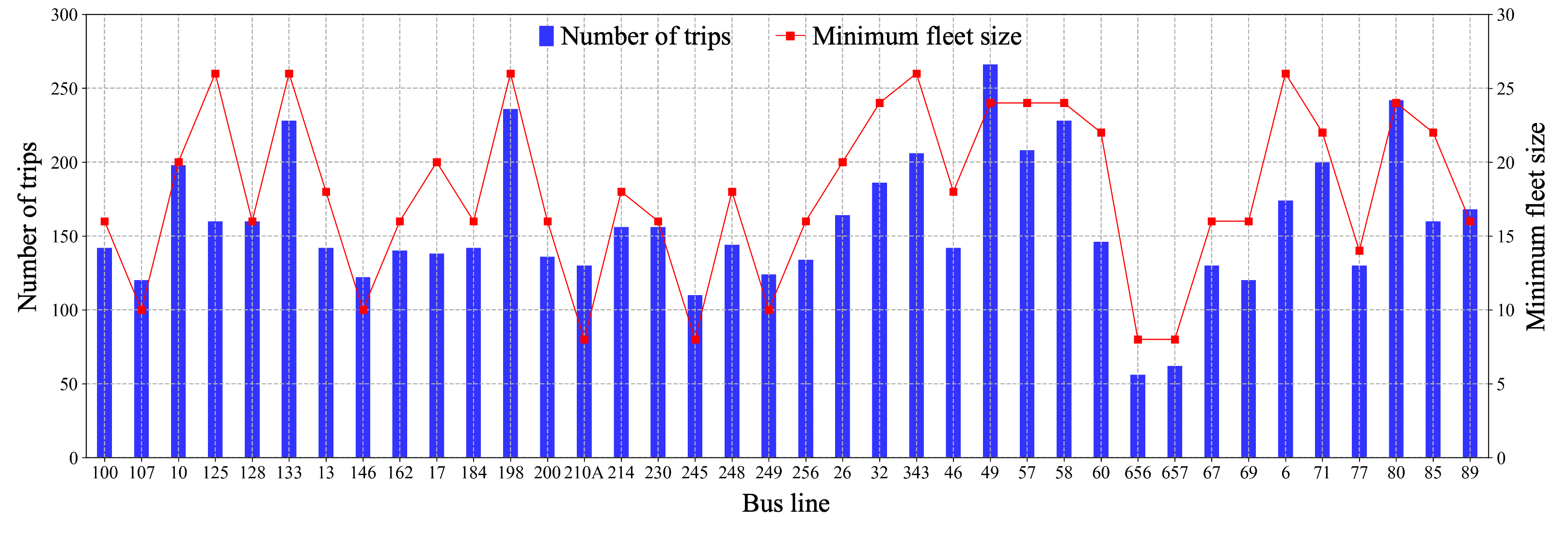}
    \caption{Number of scheduled trips (during 7:00-22:00) and minimum fleet size in each bus line.}
    \label{figtrips}
    \end{figure}

     Within the time horizon 7:00-22:00, the selected 38 bus lines have a total of 6,006 service trips, as shown in Figure \ref{figtrips}. In addition, the minimum fleet size per bus line is also plotted in this figure, following the sub-problem \eqref{FMFS1}-\eqref{FMFS6}.  
        
    \subsection{Analysis on the computational performance}
    
    All the computational performances reported below are based on a Microsoft Windows 10 platform with Intel Core i9 - 3.60GHz and 16 GB RAM. 
    
    \subsubsection{Effect of bus line pre-selection}\label{subsubsecEffls}
   As a pre-processing step,  the bus line selection module \eqref{blpseqn1}-\eqref{blpseqn3} significantly reduces the number of candidate lines. In this part, using the sequential approach, we test (1) the increase in computational efficiency, and (2) the loss of optimality, as a result of line pre-selection. The joint approach is not considered here because the lower-level problem is decoupled by line (that is, each line can be solved independently), and reducing the number of candidate lines is obviously beneficial. More importantly, the line selection module is essential to the decoupling of the lower-level problems by line, thus it is not meaningful to run the joint approach without such a preprocessing step.

   Table \ref{tabeffls2} compares the minimum fleet size module \eqref{FMFS1}-\eqref{FMFS6} with and without line selection. Both obtain global optimality while the computational time is much lower with line selection, which is expected. 
   
   \begin{table}[h!]
   \centering
   \caption{Comparison of the minimum fleet size problem with and without line pre-selection. Time refers to the computational time in seconds; Gap refers to the optimality gap provided by Gurobi.}
   \begin{tabular}{l | lllll}
   \hline
     & \# of lines & \# of trips & Min fleet size & Time (s) & Gap
     \\
     \hline
   Without & 167 & 27144 & 2824 & 1313 & 0\%
   \\
   \hline
   With & 38 & 6006 & 684 & 90 & 0\%
   \\\hline
   \end{tabular}
   \label{tabeffls2}
   \end{table}
   
   Table \ref{tabeffls} compares the sensor assignment module in the sequential approach \eqref{FOS1}-\eqref{FOS5} with and without line pre-selection, in terms of computational time and objective value $\Phi$. Note that all computational tasks have a maximum run time of 43,200s (12 hrs), beyond which the algorithm terminates and outputs the optimality gap. From this table, it is apparent that the line selection procedure significantly reduces the computational time. Regarding the loss of optimality, the objective value $\Phi$ is slightly lower with only 38 lines, as opposed to 167 lines, but only by a small margin (below 5\%).

	\begin{table}[h!]
	\centering
	\caption{Comparison of the sensor assignment module \eqref{FOS1}-\eqref{FOS5}, with and without the line pre-selection procedure. Time is the computational time in seconds, Gap refers to the optimality gap provided by Gurobi. LoO refers to the loss of Optimality calculated as the relative improvement of the objective without pre-selection. ${\bf M}$ refers to computational time reaching the maximum (43,200 s), after which the optimization procedure was terminated.}
	\begin{tabular}{c| ccc| ccc |c}
	    \hline
 \multirow{1}{*}{Total no. of}  &  \multicolumn{3}{c|}{Without pre-selection (167 lines)} & \multicolumn{4}{c}{With pre-selection (38 lines)}
 \\
 \cline{2-8}
         sensors  $N_S$     & $\Phi$ & Time (s)  & Gap & $\Phi$ & Time (s)  & Gap  & LoO
    \\\hline
    5  & 0.2325 & 107 & 0\% & 0.2229 & 1 &  0\% & 4.1\%
    \\
    10  & 0.4036 & 2808 & 0\% & 0.3838 & 10 & 0\%& 4.9\%
    \\
    15  & 0.5329 & {\bf M} & 2.71\% & 0.5168 & 325 & 0\% & 3.0\%
    \\
    20  & 0.6336 & {\bf M} & 4.50\% & 0.6193 & 4541 & 0\% & 2.3\%
\\
\hline
\end{tabular}
            \label{tabeffls} 
	\end{table}

    \subsubsection{Effect of model reduction}
    The model reduction technique in Section \ref{subsecmr} allows us to enhance the performance of the lower-level problems in the joint approach, by eliminating inefficient trip configurations. This is achieved via $\delta$, the maximum terminal idle time.  To demonstrate the effectiveness of such a technique with different values of $\delta$, Table \ref{tabCPT} presents the objective values and computational times of the 38 lower-level problems, one for each line, when the number of sensors per line is 1. By reducing $\delta$, the computational time savings are substantial (mostly within 25\%-60\%). Meanwhile, $\delta=100$ does not compromise the optimal objective values $\Phi_l$ compared to larger $\delta$, proving that large terminal idle times are indeed sub-optimal for sensing tasks.

 \setlength\LTleft{0pt}
	\setlength\LTright{0pt}
	\begin{longtable}{clcc | cl | cl | cll}
		\caption{Objective values and computational times of the lower-level problem \eqref{FAS1}-\eqref{FAS15} for each line $l\in L$ with different $\delta$ in the model reduction procedure. $\Phi_l$ is the sensing objective value for line $l$. Time is the computational time in seconds. Saving refers to the computational time reduction of $\delta=100$ relative to $\delta=300$. ${\bf M}$ refers to computational time reaching the maximum (10,000 s), after which the optimization procedure was terminated.}
            \label{tabCPT} \\
	    \hline
\multirow{3}{*}{No.}	& \multirow{2}{*}{Line} & \multirow{3}{*}{$N_{I_l}$} & \multirow{3}{*}{$N_{l}^{\text{min}}$} & \multicolumn{2}{c|}{$\delta=300$ min} & \multicolumn{2}{c|}{$\delta=200$ min} & \multicolumn{3}{c}{$\delta=100$ min}  \\
\cline{5-11}
& \multirow{2}{*}{$l\in L$} &  &  & \multirow{2}{*}{$\Phi_l$} & Time     & \multirow{2}{*}{$\Phi_l$}  & Time  & \multirow{2}{*}{$\Phi_l$} & Time & \multirow{2}{*}{Saving}
\\
&  & & &  & (s)  & &  (s) & & (s) 
\\
\hline
 1& \#100  & 142  & 16    & 0.0529  & 1335   & 0.0529 & 925  & 0.0529  & 409   & 69.4\%
\\
2& \#107  & 120  & 10  & 0.0437  & 156    & 0.0437  & 115  & 0.0437 & 93  & 40.4\%
\\
3& \#10  & 198  & 20  & 0.0461 & 5069    & 0.0461  & 4082  & 0.0461   & 1856  & 63.4\%
\\
4& \#125  & 160  & 26  & 0.0416   & 7910    & 0.0416 & 5446  & 0.0416  & 4637  & 41.4\%
\\
5& \#128  & 160  & 16  & 0.0321  & 939    & 0.0321 & 728  & 0.0321   & 520 & 44.6\%
\\ 
6& \#133  & 228  & 26  & 0.0331  & {\bf M}    & 0.0333  & {\bf M}  & 0.0346   & {\bf M} & --
\\
7& \#13  & 142  & 18  & 0.0365  & 2984    & 0.0365  & 1955  & 0.0365  & 1609 & 46.1\%
\\
8& \#146  & 122  & 10  & 0.0361  & 85    & 0.0361  & 72  & 0.0361  & 82 & 3.5\%
\\
9& \#162  & 140  & 16  & 0.0345  & 1945    & 0.0345  & 1696  & 0.0345  & 886 & 54.5\%
\\
10& \#17  & 138  & 20  & 0.0324  & 2468    & 0.0324  & 1441  & 0.0324  & 1024 & 58.5\%
\\
11& \#184  & 142  & 16  & 0.0266  & {\bf M}  & 0.0266  & {\bf M}  & 0.0266  & 3240 & --
\\
12& \#198  & 236  & 26  & 0.0345  & 7593    & 0.0345  & 4311  & 0.0345  & 3371 & 55.6\%
\\
13& \#200  & 136  & 16  & 0.0271  & 770    & 0.0271  & 643  & 0.0271  & 405 & 47.4\%
\\
14& \#210A  & 130  & 8  & 0.0182  & 8    & 0.0182  & 6  & 0.0182  & 6 & 25.0\%
\\
15& \#214  & 156  & 18  & 0.0334  & 4337    & 0.0334  & 2500  & 0.0334  & 1588 & 63.4\%
\\
16& \#230  & 156  & 16  & 0.0300  & 5946    & 0.0300  & 5021  & 0.0300  & 3958 & 33.4\%
\\
17& \#245  & 110  & 8  & 0.0253  & 13    & 0.0253  & 8  & 0.0253  & 6 & 53.9\%
\\
18& \#248  & 144  & 18  & 0.0300  & 3747    & 0.0300  & 3573  & 0.0300  & 1567 & 58.2\%
\\
19& \#249  & 124  & 10  & 0.0396  & 40    & 0.0396  & 35  & 0.0396  & 28 & 30.0\%
\\
20& \#256  & 134  & 16  & 0.0380  & 1040    & 0.0380  & 868  & 0.0380  & 786 & 24.4\%
\\
21& \#26  & 164  & 20  & 0.0566  & 1646    & 0.0566  & 895  & 0.0566  & 760 & 53.8\%
\\
22& \#32  & 186  & 24  & 0.0547  & 7259    & 0.0547  & 4357  & 0.0547  & 3440 & 52.6\%
\\
23& \#343  & 206  & 26  & 0.0601 & {\bf M} & 0.0610 & {\bf M}  & 0.0615  & 9560 & --
\\
24& \#46  & 142  & 18  & 0.0402  & 1639    & 0.0402  & 1794  & 0.0402  & 1365 & 16.7\%
\\
25& \#49  & 266  & 24  & 0.0597  &  6090  & 0.0597 & 4589  & 0.0597   & 3304 & 45.8\%
\\
26& \#57  & 208  & 24  & 0.0585  & 4853 & 0.0585 & 3496  & 0.0585  & 3291 & 32.2\%
\\
27& \#58  & 228  & 24  & 0.0464  & 7433 & 0.0464 & 7319 & 0.0464 & 5433 & 26.9\%
\\
28& \#60  & 146  & 22  & 0.0397  &  2710   & 0.0397 & 2080  & 0.0397  & 1588 & 41.4\%
\\
29& \#656  & 56  & 8  & 0.0371  & 8    & 0.0371  & 7  & 0.0371  & 5 & 37.5\%
\\
30& \#657  & 62  & 8  & 0.0240  & 4    & 0.0240  & 4  & 0.0240  & 3 & 25.0\%
\\
31& \#67  & 130  & 16  & 0.0395  & 760    & 0.0395 &  610 & 0.0395  & 474 & 37.6\%
\\
32& \#69  & 120  & 16  & 0.0381  & 1052    & 0.0381  & 658   & 0.0381  & 405  & 61.5\%
\\
33& \#6  & 174  & 26  & 0.0496 &  7686   & 0.0496 & 4438 & 0.0496  & 2881 & 62.5\%
\\
34& \#71  & 200  & 22  & 0.0371  &  {\bf M}   & 0.0371 &  {\bf M} & 0.0374  & {\bf M} & --
\\
35& \#77  & 130  & 14  & 0.0455  & 1070    & 0.0455  & 686  & 0.0455  & 678 & 36.6\%
\\
36& \#80  & 242  & 24  & 0.0506  & {\bf M}    & 0.0498  & {\bf M}  & 0.0511  & {\bf M} & --
\\
37& \#85  & 160  & 22  & 0.0372  & 3752    & 0.0372  & 2888  & 0.0372  & 2395 & 36.2\%
\\
38& \#89  & 168  & 16  & 0.0428  & 859    & 0.0428  & 648  & 0.0428  & 395 & 54.0\%
\\
\hline

\end{longtable}

    \subsection{Analysis of model solutions}
    
    \subsubsection{Overall performance of the sequential and joint approaches}
\begin{figure}[h!]
        \centering
        \includegraphics[width=\textwidth]{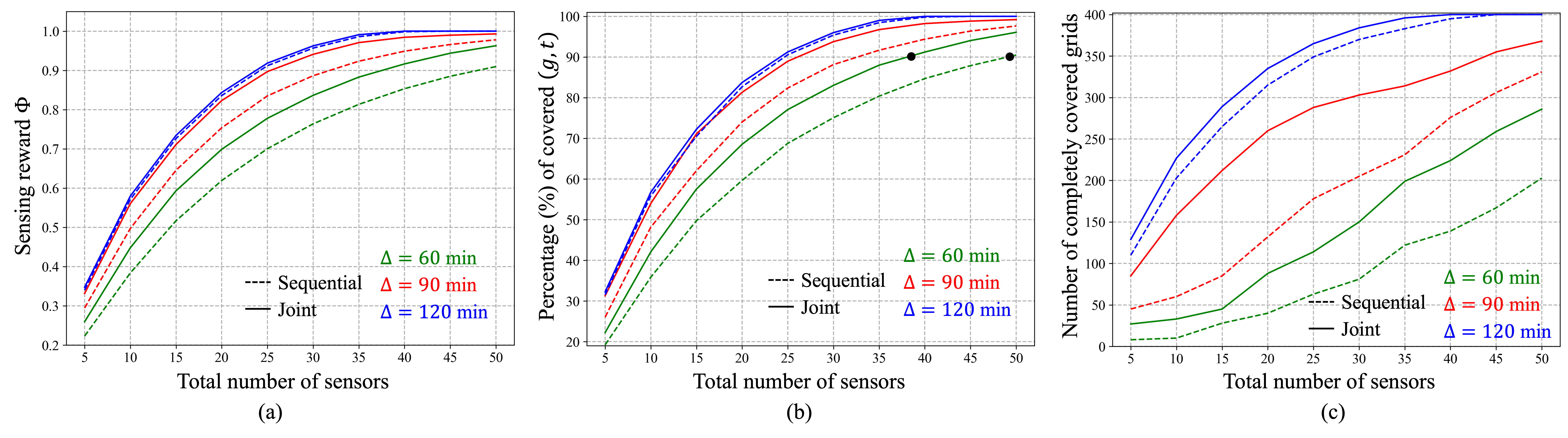}
        \caption{Sensing performance of sequential and joint approaches, in terms of sensing reward $\Phi$ (a), percentage of covered $(g,t)$ pairs (b), and number of completely covered grids (c). A completely covered grid is covered at least once in every sensing interval from 7:00-22:00.}
        \label{figcurves}
        \end{figure}

    We begin by plotting the optimal sensing rewards $\Phi$ produced by the sequential and joint approaches with given number of sensors that range from 5 to 50. The results are shown in Figure \ref{figcurves} for three different sensing intervals $\Delta=60, 90, 120$ min. It can be seen from Figure \ref{figcurves}(a) that the joint approach always yields higher $\Phi$ than the sequential approach, because the former optimizes $\Phi$ with the minimum fleet as a constraint. However, their gaps shrink as $T$ increases because higher $T$ implies less frequent sensor coverage. This suggests that the aforementioned advantage mainly lies in temporal scheduling of the bus service rather than spatial allocation of the sensors. Figure \ref{figcurves}(b) shows the percentage of covered pairs $(g,t)$, which displays a similar pattern as (a). We note that 49 sensors in the sequential approach can cover 90\% of all the target sensing subjects, while only 38 sensors are required in the joint approach. 
    
    If a spatial grid is covered in every sensing interval during 7:00-22:00, we call it {\it completely covered}. Figure \ref{figcurves}(c) compares the number of completely covered grids in the study area. It can be seen that for $\Delta=60, 90$ min, the improvement of the joint approach is substantial. In particular, for $\Delta=60$ min the increase of completely covered grids is between 41\% and 238\%. This means that coordinating bus scheduling and sensing tasks (the joint approach) yields much more reliable coverage throughout the analysis horizon.

\begin{figure}[h!]
        \centering
        \includegraphics[width=\textwidth]{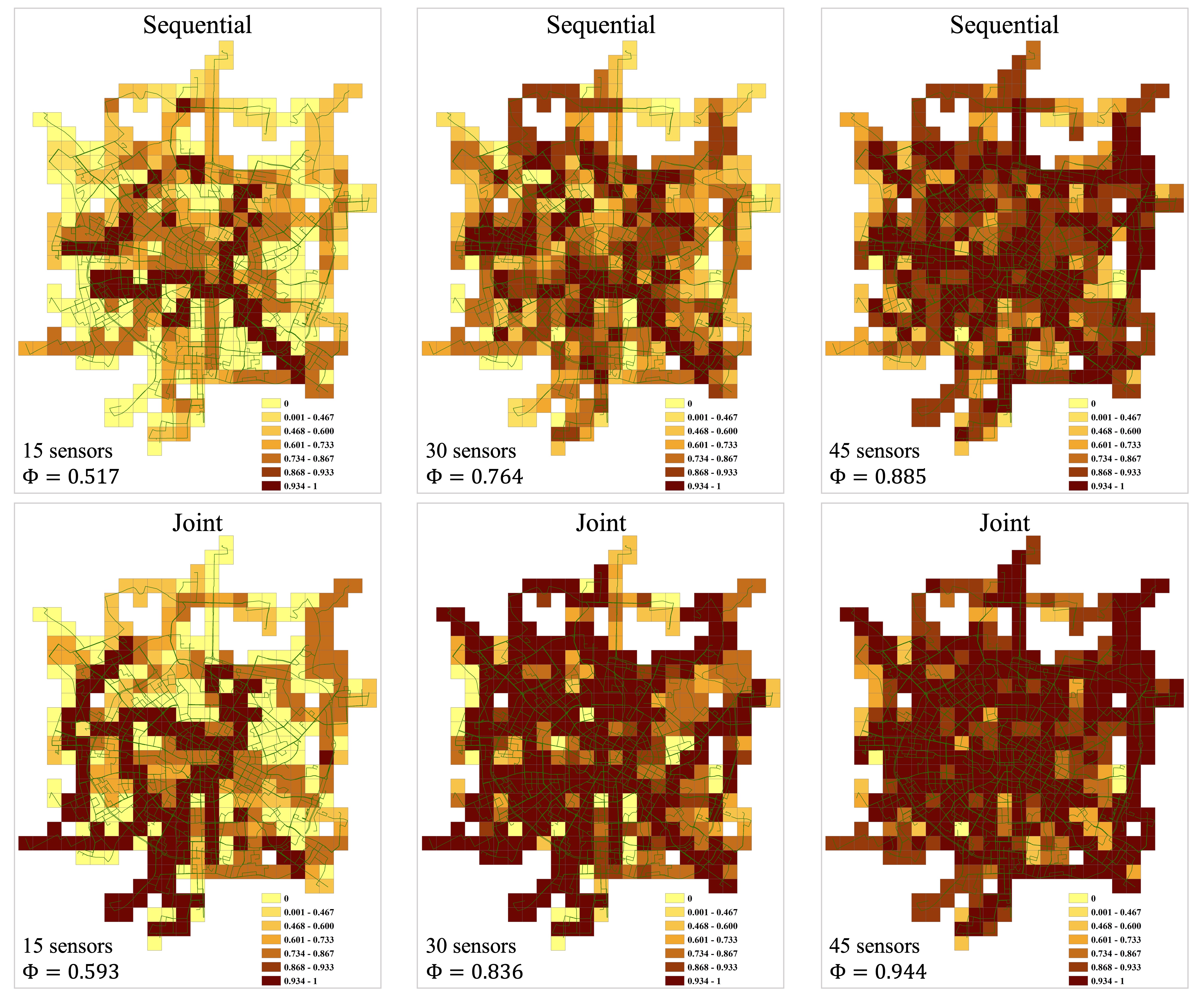}
        \caption{Visualization of the grid-based spatial coverage ($\Delta=60$ min). The value of each grid $g$ is calculated as the time-averaged grid coverage: ${1\over 15} \sum_{t=7}^{22}n_{gt}$, where $n_{gt}\in\{0, 1\}$. Grids within $[0.934, 1]$ are completely covered.}
        \label{figcoverage}
        \end{figure}

    Figure \ref{figcoverage} provides a visualization of the spatial coverage, where the grid-based values are time-averaged $n_{gt}\in\{0,\,1\}$ over 7:00-22:00, that is: ${1\over 15}\sum_{t=7}^{22} n_{gt}$. It can be seen that the coverage provided by the joint approach is much denser. 
    
    \subsubsection{Comparison of the joint approach with vehicle-based allocation} 
    As pointed out in the introduction, existing approaches on mobile sensor deployment tend to assign sensors to individual buses based on their real-world trajectories. In this part, we conduct a comparative study of the joint approach with such vehicle-based sensor allocation, in terms of optimality, computational tractability, and operational implications. 
    
\begin{figure}[H]
        \centering
        \includegraphics[width=.6\textwidth]{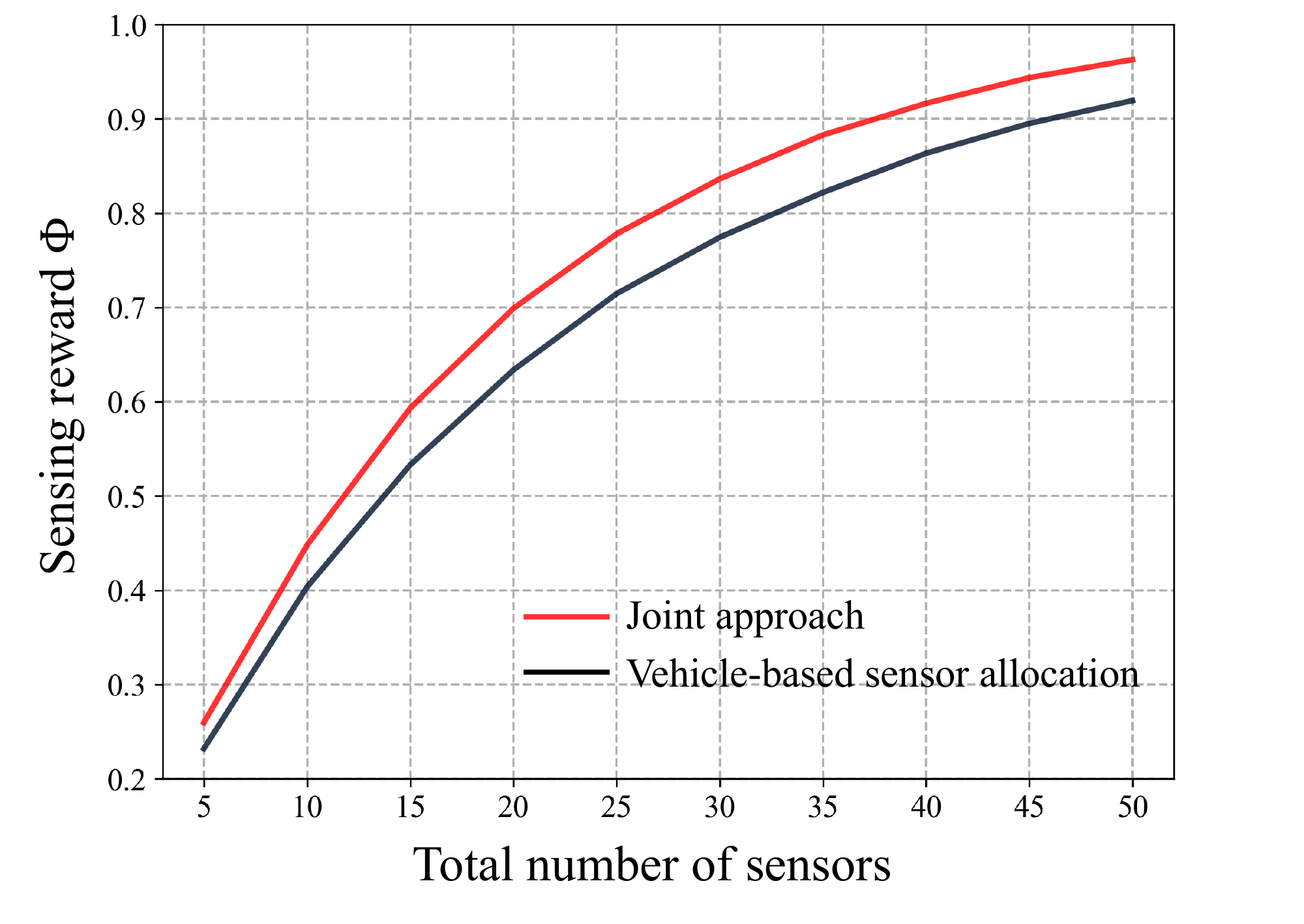}
        \caption{Solution optimality of the proposed joint approach and vehicle-based sensor allocation adopted in the literature.}
        \label{figcompare}
        \end{figure}

    For the lack of real-world bus trajectories in our case, the following scenario is used to emulate vehicle-based sensor allocation: (1) all 167 bus lines are considered without the pre-selection step (Section \ref{sec_BR_selection}); (2) the minimum-fleet-based bus scheduling (Section \ref{sec_MFS}) is performed, followed by the generation of artificial bus trajectories assuming constant speeds (this assumption has been addressed in Section \ref{subsecdiscussion}); (3) perform subset selection to the bus trajectories to maximize sensing quality. It is easy to see that this scenario is equivalent to the sequential approach without bus line pre-selection. The two approaches are compared in Figure \ref{figcompare} in terms of solution optimality, and the proposed joint approach clearly prevails.

    Regarding computation, we mainly focus on algorithm complexity as the problem size scales (i.e. with more bus lines), rather than computational time since the problem of MSD is planning in nature and not sensitive to time-efficiency of the algorithms. The following are observed: 
    \begin{itemize}
    
    \item In the joint approach, as the number of bus lines increases, the bus line pre-selection module will eliminate some redundant lines and, more importantly, the trip-combination-and-sensor allocation subproblems \eqref{FAS1}-\eqref{FAS15} are solved for each candidate line independently. Therefore, the computational burden of the joint approach is sub-linear w.r.t the number of bus lines in the problem. 
    
    \item In the vehicle-based approach, as the number of bus lines increases, the computational burden of the subset-selection problem grows combinatorically.
    
    \end{itemize}
    
    Finally, on an operational level, the vehicle-based approach does not take into account the mechanisms of bus scheduling and is prone to random noises such as traffic congestion. In contrast, the trip-based joint approach can easily adapt to operational changes or objectives, and is robust against service uncertainties per our discussion in Section \ref{subsecdiscussion}.

        \subsubsection{Non-uniqueness of the minimum fleet solution}
    
    One of the main tricks of the joint optimization approach is to treat the minimum fleet size as a constraint for each line, and finding the trip chain configuration, among all possible ones, that promises the best sensing outcome.  In this section, we analyze the non-uniqueness of the minimum fleet solution. Three sample lines (\#200, \#248, \#184) are chosen to illustrate their trip chain configurations from the sequential and joint approaches. As shown in Table \ref{tab_EAS}, their minimum fleet sizes are 16, 18, and 16, respectively. Assuming only 1 sensor per line, we compare the number of trips served and $(g,t)$ pairs covered by the instrumented trip chain in the sequential and joint solutions, and find that the joint approach outperforms the sequential approach in terms of both in all three lines. This shows that different trip chain sets can indeed differ considerably in sensing performance, even though they have the same fleet size.

\begin{figure}[h!]
    \centering
    \includegraphics[width=1.0\textwidth]{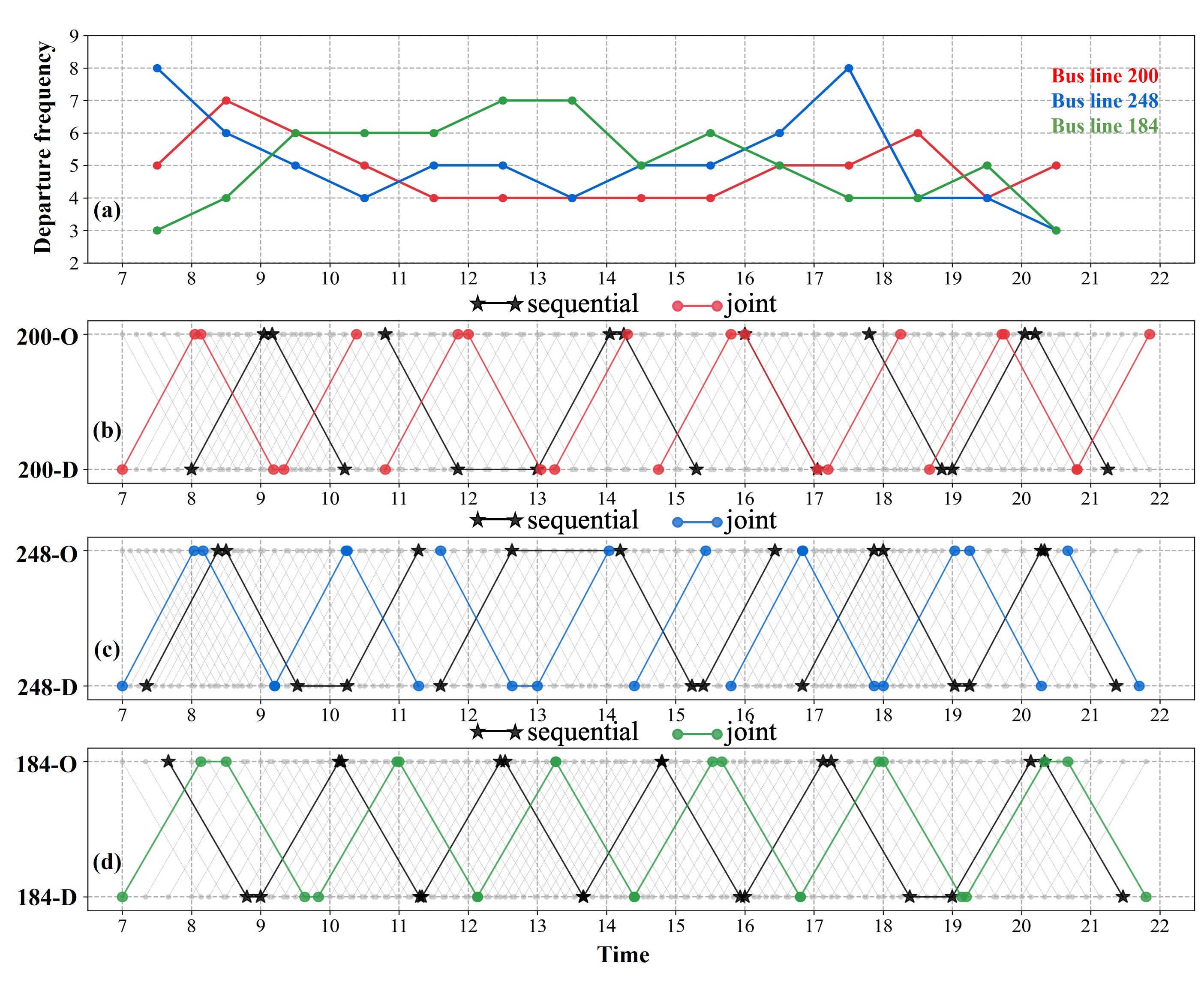}
    \caption{Visualization of instrumented trip chains in the sequential and joint approaches. (a) Departure frequencies of bus lines \#200, \#248, and \#184; (b)-(d) instrumented trip chains in line \#200, \#248, and \#184.}
    \label{fig_SL_timetable}
    \end{figure}

  These instrumented trip chains are further visualized in Figure \ref{fig_SL_timetable}. It can be seen that the joint approach indeed covers more trips, with very little terminal idling. On the other hand, in the sequential approach the instrumented buses were idle during 11:50-13:00 (\#200), 12:40-14:10 (\#248) and 18:20-19:00 (\#184), which are all within the off-peaks of the respective bus lines; see Figure \ref{fig_SL_timetable}(a). This suggests that without coordinating trip chain formation (i.e. bus scheduling) and sensing tasks, the sequential approach tends to under-utilize the sensing capacity of instrumented buses.

    \begin{table}[h!]
    \centering
    \caption{The effect of different minimum fleet solutions on the sensing outcome.}
    \begin{tabular}{cccccc}
    \hline
    \multirow{2}{*}{Bus line} & Solution & \multirow{2}{*}{Min fleet size} & Number of & Number of  & Number of $(g,t)$ 
    \\
    & procedure & &  sensors & trips served & pairs covered
     \\ 
     \hline
   \multirow{2}{*}{\#200} & Sequential   & \multirow{2}{*}{16}   & \multirow{2}{*}{1}     & 9     & 125      
     \\
                                     & Joint           &       &     & 12     & 160                                                                                            
    \\
    \hline
   \multirow{2}{*}{\#248}    & Sequential   & \multirow{2}{*}{18}   & \multirow{2}{*}{1}   & 10     & 144                                                                                           \\
                                          & Joint          &     &      & 12     & 179                                                                                            
    \\
    \hline
    \multirow{2}{*}{\#184}     & Sequential   & \multirow{2}{*}{16}    &  \multirow{2}{*}{1}   & 11   & 136                                                                                             \\
                                            & Joint           &      &     & 12  & 145                                                                                            
        \\ \hline
    \end{tabular}
    \label{tab_EAS}
    \end{table}

%

\begin{figure}[h!]
    \centering
    \includegraphics[width=\textwidth]{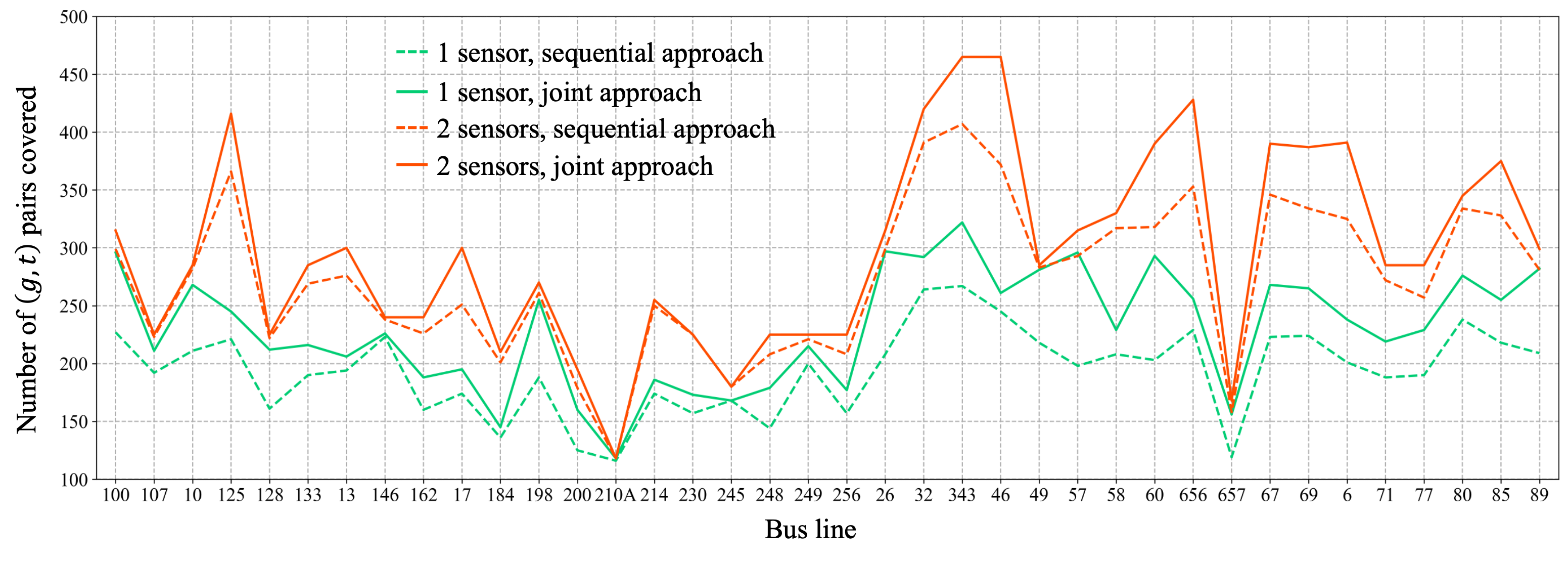}
    \caption{Comparison of sequential and joint approaches in terms of number of $(g,t)$ pairs covered per line.}
    \label{figpairscov}
    \end{figure}

Finally, Figure \ref{figpairscov} provides an overview of the sensing enhancement by considering the non-uniqueness of the minimum fleet solution. The joint approach provides much more coverage than the sequential one for most of the bus lines, and their gaps are generally larger for 1 sensor than 2 sensors.

    \subsubsection{Line saturation}\label{subsubsectionsaturate}

In this section, we investigate the number of sensors required to saturate a line in the joint trip-chain-configuration-and-sensor-allocation problem. A line is called saturated with $K$ sensors if $K$ is the smallest integer such that the sensing objective with $K+1$ sensors does not improve. This issue is of interest to us because the joint approach requires repeatedly solving the lower-level problem \eqref{FAS1}-\eqref{FAS15} for $m_l=0, 1, \ldots, K_l$; see \eqref{FSTR5'}. It is critical to find the upper bound $K_l$ such that line $l$ is saturated.

\begin{figure}[h!]
    \centering
    \includegraphics[width=1\textwidth]{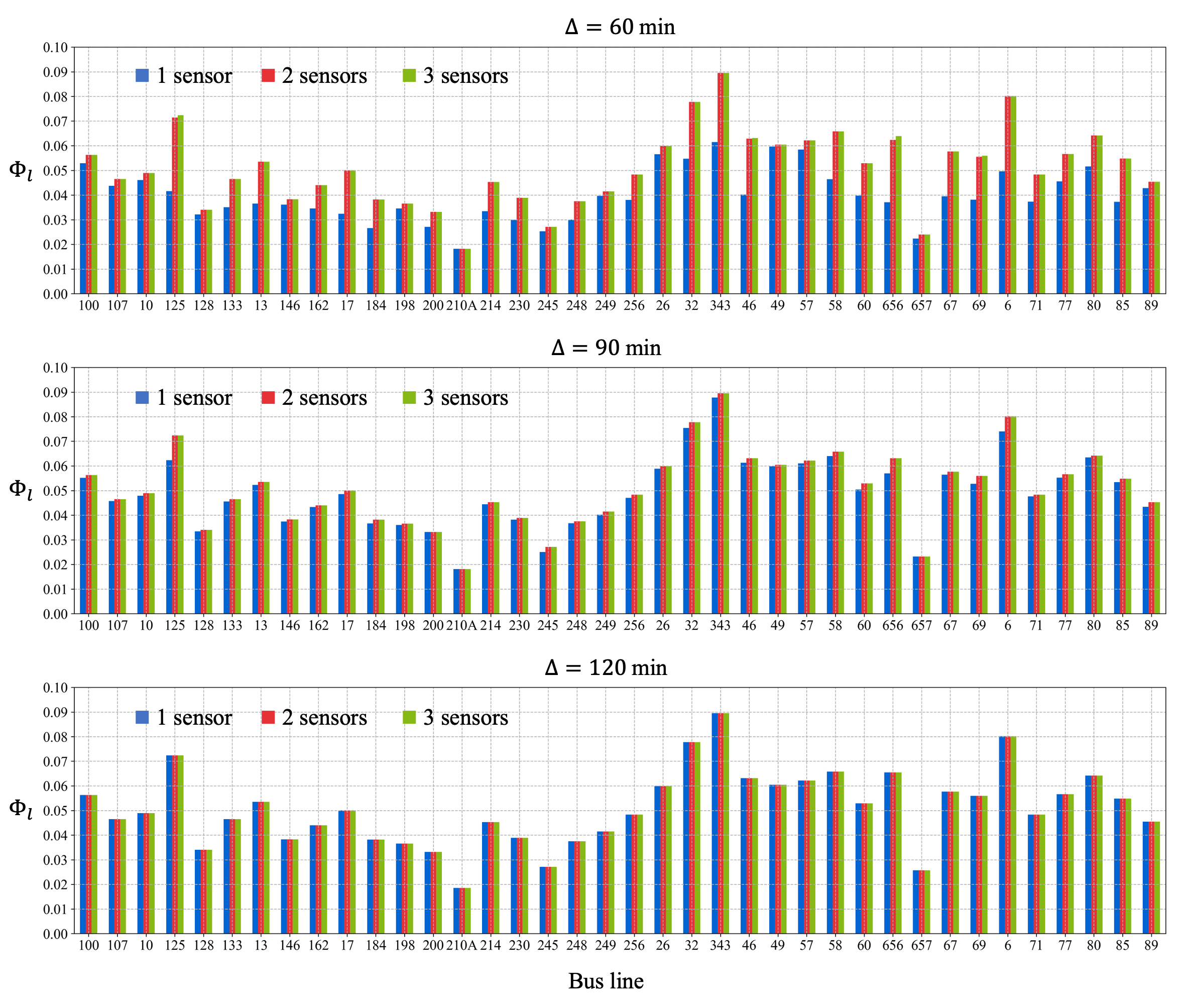}
    \caption{Sensing objective $\Phi_l$ for each of the 38 lines when $1$, $2$ and $3$ sensors are respectively installed per line. Each sub-figure corresponds to a sensing interval length $\Delta$.}
    \label{figsaturation}
    \end{figure}

Such bounds are inferred in Figure \ref{figsaturation}: 2 sensors are sufficient to saturate most of the lines when $\Delta=60$ min, and all the lines when $\Delta=90$ min. For $\Delta=120$ min, 1 sensor is sufficient. Therefore, a general rule of thumb is that 2 sensors per line should be adequate for the type of sensing objective introduced in this paper. This also confirms that the computational burden of the upper-level problem is insignificant.

            \begin{figure}[h!]
    \centering
    \includegraphics[width=0.55\textwidth]{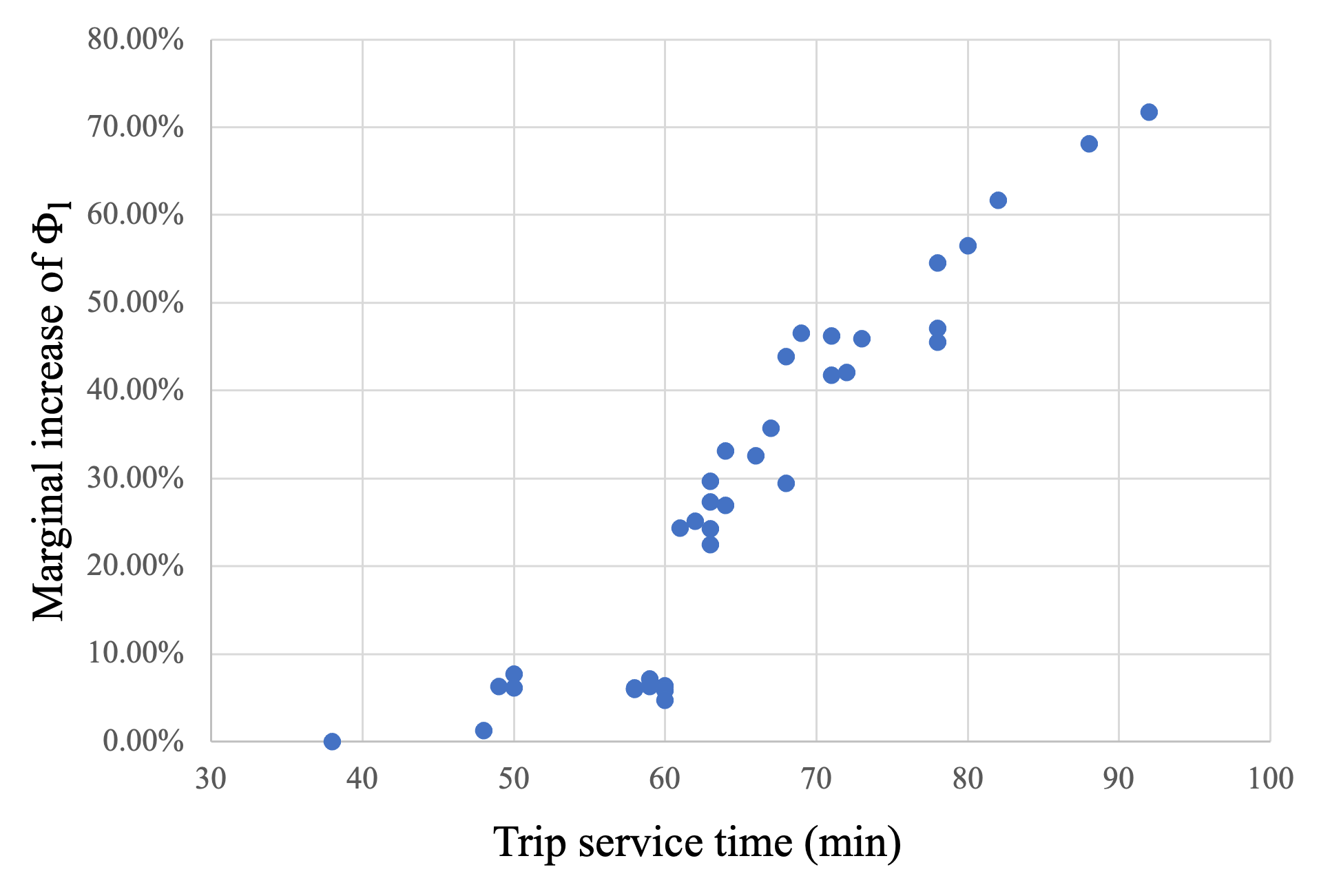}
    \caption{Scatter plot of trip service time against the marginal gain of $\Phi_l$ by having 2 sensors instead of 1.}
    \label{figgain}
    \end{figure}

    For $\Delta=60$ min, some bus lines benefit greatly from having a second sensor, while some do not. Figure \ref{figgain} plots such marginal gains of $\Phi_l$ (from having 1 to 2 sensors) against average trip service times (one-way run time) of the 38 lines.  A clear monotonic relationship exists between the two, suggesting that lines with large run times usually require a second sensor to be saturated. This is understandable given that the average time gap between two consecutive visits by the same sensor experienced by a grid is equal to the trip service time; the longer the service time, the more grids left uncovered in an hour there will be. Such a finding suggests a way to quickly generate approximate solutions to the MSD problem; see Section \ref{sec_conclusion}.

\section{Conclusion and discussion} \label{sec_conclusion}

This paper addresses the mobile sensor deployment problem in bus transit systems. A trip-based sensor deployment framework is proposed, which effectively integrates bus scheduling with mobile sensing. To solve this problem in large-scale instances, we propose a multi-stage optimization procedure, which includes line selection, minimum fleet size, and sensor allocation sub-modules. To further explore the coordination between bus scheduling and sensing tasks, we propose a bi-level formulation for the joint optimization of fleet size and sensing quality. Specific findings and discussions are as follows. 

\begin{itemize}

\item The line selection module, as a preprocessing step, can considerably reduce the computational burden in subsequent procedures. It is also essential to the joint approach by decoupling the simultaneous trip-chain-configuration-and-sensor-allocation problems by each line, which can be solved in a distributed fashion.

\item  The model reduction technique in Section \ref{subsecmr} allows us to enhance the performance of the lower-level problems in the joint approach, by eliminating inefficient trip configurations. Specifically, this is done by imposing the artificial constraint \eqref{eqnarti} and eliminating the big-M constraint \eqref{FAS2}. As a result, the computational time savings in the lower-level problems are mostly between 25\%-60\%.

\item The joint approach takes advantage of the non-uniqueness of the minimum fleet solution, by coordinating trip chain formation (i.e. bus scheduling) and sensing tasks. This results in considerable sensing enhancement in the case study: 49 sensors can cover 90\% of the sensing subjects $(g,t)$ in the sequential approach, and only 38 sensors are sufficient for this in the joint approach. In addition, the number of grids covered in every hour throughout the time horizon has increased by 41\%-238\% in the joint approach. 

\item The joint approach prevails because it minimizes terminal idling times for the instrumented buses, whereas the sequential approach suffers from long idling times -- a consequence of pursuing minimum fleet size (which corresponds to multiple trip chain configurations) before considering sensing efficacy. Such a difference is pronounced for high-frequency sensing scenarios (e.g. air quality, traffic conditions). However, the advantage of the joint approach diminishes as $\Delta$ increases (e.g. heat island phenomena). For very large $\Delta$ (e.g. road surface condition), the temporal aspect of the optimization can be ignored and the MSD problem reduces to assigning sensors to bus lines. A recent study \citep{DH2023} considers inter-line relocations of instrumented buses, which are suitable for this kind of applications.

\item For all the bus lines in the case study, it is found that 2 sensors are sufficient to saturate a line when $\Delta=60$ min, in the sense that a third one brings very limited gain in sensing rewards. For longer sensing intervals $\Delta=90, 120$ min, 1 sensor is sufficient. This suggests a way to quickly generate approximate solutions (i.e. prototyping) without being computationally involved: Taking $\Delta=60$ min as an example, apply line selection followed by allocating 1 sensor to each, and extra ones can be allocated to those with long trip service times. Then, bus scheduling and sensor allocation can be performed within each line independently, and even a little deviation from our theoretical solutions (for various practical reasons) would have limited impact on the overall sensing quality.
\end{itemize}

\section*{Acknowledgement}
This work is supported by the National Natural Science Foundation of China through grants 72071163 and 72271206, and the Natural Science Foundation of Sichuan Province through grant 2022NSFSC0474.

\end{document}